\documentclass[onefignum,onetabnum]{siamart171218}

\usepackage[T2A]{fontenc}
\usepackage[cp866]{inputenc}
\usepackage[english]{babel}

\usepackage{amssymb}
\usepackage{tikz-cd}
\usepackage{amsmath}
\usepackage{latexsym}
\usepackage{amsfonts}
\usepackage{graphicx}
\usepackage{wrapfig}
\usepackage{bbm}

\DeclareMathOperator{\diag}{diag}

\DeclareMathOperator{\sech}{sech}

\DeclareMathOperator{\sgn}{sgn}
\DeclareMathOperator{\supp}{supp}
\DeclareMathOperator{\dist}{dist}

\DeclareMathOperator{\Rea}{Re}
\DeclareMathOperator{\Ima}{Im}

\newcommand{\loc}{\text{loc}}

\newsiamremark{remark}{Remark}

\numberwithin{equation}{section}

\headers{}{S. Shindin, N. Parumasur, and O. Aluko}

\title{Analysis of Malmquist-Takenaka-Christov rational approximations
with applications to the nonlinear Benjamin equation}

\author{Sergey Shindin\thanks{School of Mathematics, Statistics and Computer Sciences,
University of KwaZulu-Natal (\email{shindins@ukzn.ac.za}).}
\and Nabendra Parumasur\thanks{School of Mathematics, Statistics and Computer Sciences,
University of KwaZulu-Natal (\email{parumasurn1@ukzn.ac.za}).}
\and Olabisi Aluko\thanks{School of Mathematics, Statistics and Computer Sciences,
University of KwaZulu-Natal (\email{olabisialuko9.oa@gmail.com}).}}

\begin{document}
\maketitle

\begin{abstract}
In the paper, we study approximation properties of the Malmquist-Takenaka-Christ\-ov (MTC) 
system. We show that the discrete MTC approximations converge rapidly under mild restrictions 
on functions asymptotic at infinity. This makes them particularly suitable for solving semi- and 
quasi-linear problems containing Fourier multipliers, whose symbols are not smooth at the origin.
As a concrete application, we provide rigorous convergence and stability analyses of 
a collocation MTC scheme for solving the nonlinear Benjamin equation. 
We demonstrate that the method converges rapidly and admits an efficient implementation, 
comparable to the best spectral Fourier and hybrid spectral Fourier/finite-element methods 
described in the literature.
\end{abstract}



\section{Introduction}\label{sec1}
In the paper, we consider 
the nonlinear equation, proposed by T.B. Benjamin in his study of internal waves arising in a two fluid system, 
see \cite{BENJ1996}. The equation reads
\begin{subequations}\label{eq1.1}
\begin{equation}\label{eq1.1a}
u_t =-\alpha u_x+\beta\mathcal{H}[u_{xx}]+\gamma u_{xxx} -\delta (u^2)_x, \quad u(0)=u_0, 
\quad x\in\mathbb{R},
\end{equation}
where $\alpha$, $\beta$, $\gamma$ and $\delta$ are real parameters and
\[
\mathcal{H}[u](x) = \tfrac{\text{p.v.}}{\pi}\int_{\mathbb{R}} \frac{u(y)}{y-x} dy,
\]
is the standard Hilbert transform. The problem is formally Hamiltonian, i.e.
\begin{equation}\label{eq1.1b}
u_t = \mathcal{J} \nabla \mathcal{G}(u),
\end{equation}
where $\mathcal{J} = - \partial_x$ is the skew-symmetric first order automorphism of the Hilbert scale 
$H^{s}(\mathbb{R})$, $s\in\mathbb{R}$ and 
\begin{equation}\label{eq1.1c}
\mathcal{G}(u)=\tfrac{1}{2}\int_{\mathbb{R}}
\Bigl(\alpha|u|^2 - \beta u\mathcal{H}[u_x]+\gamma |u_{x}|^2  +  \tfrac{2\delta}{3}u^{3}\Bigr) dx.
\end{equation}
\end{subequations}

In recent years, problem \eqref{eq1.1} has received significant attention in both analytic and numerical 
communities. The well-posedness analysis of \eqref{eq1.1a} can be found in 
\cite{ChGuXi2011, LIWU2010, LIN1999}. In particular, the arguments of 
\cite{ChGuXi2011}, indicate that \eqref{eq1.1a} is globally well-posed, provided the initial data $u_0$ 
is in $H^{s}(\mathbb{R})$, with $s\ge -\frac{3}{4}$. The global classical solutions are obtained if 
$u_0\in H^{s}(\mathbb{R})$ and $s\ge 3$. The study of traveling wave solutions is initiated 
in \cite{BENJ1996}. The existence of such solutions for all admissible values of the model parameters 
is affirmatively settled by several authors (see e.g. \cite{ChBo1998, Pava1999} and references therein), 
while their stability is discussed in \cite{AlBoRe1999, BENJ1996, Pava1999}.

On the numerical side, a variety of techniques, suitable for integrating \eqref{eq1.1}, as well as 
for finding associated traveling waves, is described in the literature. Among others, we mention 
the pseudo-spectral Fourier-type schemes used in \cite{AlBoRe1999, KaBo2000}, the hybrid 
Fourier-type/finite-difference scheme of \cite{CaAk2003} and the hybrid Fourier-type/finite-element 
methods employed in \cite{DoDuMi2015, DoDuMi2016}.
In all the techniques listed above, the spatial domain is truncated to a large interval $[-L,L]$, the 
resulting stationary and/or non-stationary Benjamin equation, equipped with periodic boundary conditions, 
is solved numerically. However, as observed by a number of authors, due to the jump discontinuity in the Fourier 
symbol of the operator $\mathcal{H}$, the exact solutions decay at most algebraically at 
infinity.\footnote{For a recent study of 
the interplay between regularity and asymptotic of solutions see \cite{Urrea2013}.} This is a serious technical 
obstacle as an accurate numerical approximation of such solutions requires very large values for the 
truncation parameter $L$ (see e.g. the discussion and numerical experiments in \cite{CaAk2003, DoDuMi2015}).

In the paper, we adopt an alternative approach. We approximate solutions directly in the real line using 
a family of rational orthogonal functions proposed independently by F. Malmquist \cite{Malm1926}, 
S. Takenaka \cite{Tak1926} and, in context of spectral methods, by C.I. Christov \cite{Cri1982}.
The Malmquist-Takenaka-Christov (MTC) system has a number of attractive computational features. 
As observed by J. Weideman \cite{Wid1995}, the MTC function are eigenfunctions of the Hilbert transform;
the system behaves well with respect to the product of its members \cite{Cri1982}; the MTC
differentiation matrices are skew symmetric and tridiagonal, while computing of the discrete spectral 
MTC coefficients can be accomplished efficiently via discrete Fast Fourier Transform (FFT) 
\cite{Cri1982, Wid1994, Wid1995}. In fact, it is observed recently in \cite{IsWe2019} that the MTC system is  
the only complete rational orthogonal basis in $L^2(\mathbb{R})$ that possesses the last two properties. 

Unfortunately, not much is known about the convergence rate of the MTC-Fourier series. Some preliminary 
results in this direction are obtained in \cite{Boyd1990, Wid1995}, where it is shown that the 
convergence rate is geometric, provided functions under consideration are analytic in an exterior 
of some neighborhood of $\{i,-i\}$ in the complex plane. However, as noted in \cite{IsWe2019, Wid1995}, 
these results have limited applications, specifically in the context of differential equations. 

In the paper, we derive several error bounds describing convergence of the continuous and discrete  
MTC-Fourier expansions. It turns out that the convergence rate is controlled solely 
by the regularity and asymptotics of the Fourier images of functions in $\mathbb{R}\setminus\{0\}$, while allowing 
square integrable singularities at the origin. As a consequence, and in contrast to the Hermite 
or algebraically mapped Chebyshev bases \cite{Boyd, Canuto, Fun1992} in $L^2(\mathbb{R})$, the 
MTC-Fourier approximations converge spectrally under very mild restrictions on the functions decay at infinity. 
The latter circumstance makes them particularly suitable when dealing with semi- and quasi-linear equations 
containing Fourier multipliers, whose symbols are not smooth at the origin (e.g. the Hilbert/Riesz 
transforms, fractional derivatives, e.t.c.). In the concrete case of the Benjamin equation \eqref{eq1.1}, 
the MTC semi-discretization  yields a spectrally convergent collocation scheme that admits an efficient 
practical implementation, comparable to the best spectral-Fourier and hybrid spectral-Fourier/finite-element 
methods, described in the literature.\footnote{Similar technique is employed recently in 
\cite{BoXu2011, BoXu2012} for the closely related Benjamin-Ono equation.} 

The paper is organized as follows: In Section~\ref{sec2}, we fix the notation and provide a basic function theoretic
setup that is used in our analysis. Section~\ref{sec2} contains a technical result, for which we have no immediate
references. For the readers convenience, we sketched the proof in Appendix~\ref{sec7}.
A detailed discussion of the MTC basis and its approximation properties 
is the main subject of Section~\ref{sec3}. A convergence analysis of the MTC collocation scheme 
suitable for numerical integration of \eqref{eq1.1} is presented in Section~\ref{sec4}. 
Numerical experiments, illustrating computational performance of our scheme, are reported in 
Section~\ref{sec5}. Section~\ref{sec6} is reserved for concluding remarks.

\section{Preliminaries}\label{sec2}

This section is introductory. Here, we fix a notation and provide a basic function theoretic setup pertinent for  
our calculations.

\paragraph{Notation} Throughout the paper, symbols
\begin{align*}
&\mathcal{F}[\varphi](\xi) = \hat{\varphi}(\xi)= \tfrac{1}{\sqrt{2\pi}}\int_{\mathbb{R}} e^{-i\xi x}\varphi(x)dx,\\
&\mathcal{F}^{-1}[\hat{\varphi}](x) = \varphi(x)= \tfrac{1}{\sqrt{2\pi}}\int_{\mathbb{R}} e^{i\xi x}
\hat{\varphi}(\xi)d\xi,
\end{align*}
denote the normalized Fourier transforms and its inverse. Letters $x$ and $\xi$ are reserved for 
the physical and the frequency variables, respectively. Symbol $\ast$ denotes the standard 
Fourier convolution. Letter $c$ stands for a generic positive constants, whose particular value is irrelevant.

\paragraph{Weighted Lebesgue spaces} 
Let $\Omega\subseteq \mathbb{R}$ be measurable and let $w\in L^{1,\loc}(\Omega)$ 
be a.e. positive in $\Omega$. We employ
\[
L^p_w(\Omega; B) := L^p(\Omega, wdx; B),\quad 1\le p\le \infty,
\]
to denote weighted Lebesgue spaces with values in a Banach space $B$, we write shortly $L^p_w(\Omega)$, 
when $B $ is either of $\mathbb{R}$ or $\mathbb{C}$. In the sequel, we deal with power weights 
$w_{pr}(x) = x^{pr}$, $r>-\frac{1}{p}$. For such weights, we use the shortcut 
$L^p_r(\mathbb{R}_+)$. When $r=0$, we write simply $L^p(\Omega)$.

\paragraph{Variable weight Sobolev spaces} 
The error analysis of Section~\ref{sec3} in a natural way gives rise to a scale of variable weight Sobolev spaces.  
For real valued functions, these are defined by\footnote{For $s\in\mathbb{N}$, 
the $H^s_r(\mathbb{R})$-norm is equivalent to $\|f\|_{L^2(\mathbb{R})}+\sum_{m=0}^s 
\|\kappa^{\pm}_{r-s+m,\ell} f^{(m)}\|_{L^2(\mathbb{R})}$, i.e. $\|\cdot\|_{H^s_r(\mathbb{R})}$ 
is a Sobolev-like norm, where weak derivatives of different orders are integrated against different weights, 
hence the name.} 
\begin{subequations}\label{eq2.1}
\begin{align}
\label{eq2.1a}
H^s_r(\mathbb{R}) 
&= \{f\in \Rea \mathcal{S}' \,|\, \|f\|_{H^s_r(\mathbb{R})}<\infty\},
\quad s>-\tfrac{1}{2},\quad r\ge 0,\\
\label{eq2.1b}
\|f\|_{H^s_r(\mathbb{R})}^2
&= \tfrac{1}{2}\|f\|_{L^2(\mathbb{R})}^2
+\|\mathcal{P}_\pm[\kappa_{r,\ell}^\pm f]\|_{\mathring{H}^{s}(\mathbb{R})}^2,
\end{align}
\end{subequations}
where $\Rea \mathcal{S}'$ is the space of real valued tempered distributions, $\mathcal{P}_{\pm}$ 
are Fourier multiplier (projectors) associated with the Heaviside functions 
$\hat{h}^\pm(\xi) = \tfrac{1\pm \sgn(\xi)}{2}$, $\kappa_{r,\ell}^\pm(\cdot) 
= \tfrac{1}{\sqrt{2\pi}}(i\ell \pm \cdot)^r$ and $\mathring{H}^{s}(\mathbb{R})$ 
is the standard homogeneous Sobolev space of order $s$, see e.g. \cite{BeLo1976}. 
The meaning of parameters $s$, $r$ and $\ell$ is straightforward.  Parameter $s>0$ controls 
regularity of $f$, while $r$ describes its asymptotic at infinity. The positive scaling parameter 
$\ell$ is used in practical simulations to control the distribution of spatial nodes and to tune 
up the convergence rate. 

Basic properties of the variable weight Sobolev space $H^s_r(\mathbb{R})$ are contained in the following 
\begin{lemma}\label{lm2.1}
$H^s_r(\mathbb{R})$, with $s>-\tfrac{1}{2}$ and $r\ge 0$, are Hilbert spaces. Further,
\begin{subequations}\label{eq2.2}
\begin{equation}\label{eq2.2a}
H^{s_0}_{r_0}(\mathbb{R}) \hookrightarrow H^{s_1}_{r_1}(\mathbb{R}),
\end{equation}
provided  
\begin{equation}\label{eq2.2b}
-\frac{1}{2}<s_1\leq s_0 \leq s_1+r_0-r_1,\quad  0\le r_1\le r_0.
\end{equation}
\end{subequations}
Finally, for $s_0,s_1>-\tfrac{1}{2}$, $r_0,r_1\ge 0$ and $\theta\in(0,1)$, we have
\begin{equation}\label{eq2.3}
[H^{s_0}_{r_0}(\mathbb{R}), H^{s_1}_{r_1}(\mathbb{R})]_\theta 
= H^{(1-\theta)s_0+\theta s_1}_{(1-\theta)r_0+\theta r_1}(\mathbb{R}),
\end{equation}
where $[\cdot,\cdot]_\theta$ denotes the standard complex interpolation functor of A. Calderon \cite{BeLo1976}.
\end{lemma}
\begin{proof}
(a) In terms of Fourier images, \eqref{eq2.1b} reads
\begin{equation}\label{eq2.4}
\|f\|_{H^s_r(\mathbb{R})}^2 = \|\hat{f}\|_{L^2(\mathbb{R}_\pm)}^2
+\|\mathcal{J}_{\mp}^{-r}[\hat{f}]\|_{L^2_s(\mathbb{R}_\pm)}^2,
\quad \mathcal{J}_{\mp}^{-r}[\hat{f}] = \sqrt{2\pi}\bigl(\hat{\kappa}^\pm_{r,\ell}\ast \hat{f} \bigr).
\end{equation}
Note that $\supp\hat{\kappa}^{\pm}_{r,\ell}\subset\mathbb{R}_{\mp}$. Hence, for real 
valued distributions (whose Fourier images are Hermitian) the choice of sign in \eqref{eq2.1b}, \eqref{eq2.4} 
is irrelevant. 

Operators $\mathcal{J}_{\mp}^{r}[\cdot](\xi)$, $\xi\in \mathbb{R}_\pm$, defined in \eqref{eq2.4}, 
are known as one-sided Bessel potentials in the half line, \cite{Samko}. In the context of the Laguerre basis in 
$\mathbb{R}_+$, weighted spaces of such potentials are discussed in \cite{BaPaPoSh2014}. 
In particular, it is shown that 
\[
L^{2,r}_s(\mathbb{R}_\pm) :=\mathcal{J}_{\mp}^r[L^2_s(\mathbb{R}_\pm)], 
\quad s>-\tfrac{1}{2},\quad r\ge 0, 
\]
equipped with the norm 
\[
\|\cdot\|_{L^{2,r}_s(\mathbb{R}_\pm)} = \|\mathcal{J}_\mp^{-r}[\cdot]\|_{L^2_s(\mathbb{R}_\pm)}, 
\]
are Banach spaces.\footnote{ In fact, only the case $\ell=1$ and $\mathbb{R}_+$ is treated their, 
but the extension to $\mathbb{R}_-$ and arbitrary $\ell>0$ is straightforward.}
Since $L^2(\mathbb{R})\cap \Rea \mathcal{S}'$ distributions are regular, we conclude that the quantity 
$\|\cdot \|_{H^s_r(\mathbb{R})}$ is a norm in $H^s_r(\mathbb{R})$.
The completeness of $H^s_r(\mathbb{R})$ follows from the completeness of 
$L^{2,r}_s(\mathbb{R}_\pm)\cap L^2(\mathbb{R}_\pm)$. In view of \eqref{eq2.4}, the bilinear form
\[
\langle f, g \rangle_{H^s_r(\mathbb{R})} 
= \langle \hat{f},\hat{g} \rangle_{L^2(\mathbb{R}_+)}
+\langle \mathcal{J}^{-r}_{\mp}[\hat{f}],
\mathcal{J}^{-r}_{\mp}[\hat{g}] \rangle_{L^2_s(\mathbb{R}_+)}
\]
is the inner product in $H^s_r(\mathbb{R})$. Hence, the first claim of Lemma~\ref{lm2.1} is settled.

(b) Embedding \eqref{eq2.2} is the direct consequence of \eqref{eq2.4} and the embedding inequality (21) 
from \cite{BaPaPoSh2014}.

(c) Interpolation identity \eqref{eq2.3} follow from \cite[Theorem 5.6.3]{BeLo1976} and formula
\begin{equation}\label{eq2.5}
[L^{p,r_0}_{s_0}(\mathbb{R}_+),L^{p,r_1}_{s_1}(\mathbb{R}_+)]_\theta 
= L^{p,(1-\theta)r_0 + \theta r_1}_{(1-\theta)s_0+\theta s_1}(\mathbb{R}_+),
\quad 1<p<\infty,
\end{equation}
if we view $H^s_r(\mathbb{R})$ as a retract of the vector valued Banach space 
$\tilde{H}^s_r(\mathbb{R})=\{(u,v)| u\in L^2(\mathbb{R}),\, \hat{v} \in L^{2,r}_s(\mathbb{R}_-)\cap
L^{2,r}_s(\mathbb{R}_+)\}$.

We remark that for $s_0=s_1=0$, \eqref{eq2.5} is well known. 
Proving \eqref{eq2.5} in general, even in the basic settings of the half-line, is a delicate issue\footnote{
Specifically, when either of $s_0$, $s_1$ is outside the range $\bigl(-\tfrac{1}{p},\tfrac{1}{p'}\bigr)$.}
which is of an independent interest. For readers convenience, we sketch the proof in Appendix~\ref{sec7}. 
\end{proof}

To conclude this section, we note that $H^{s}_0(\mathbb{R})=H^s(\mathbb{R})$, 
where $H^s(\mathbb{R})$ is the standard Sobolev spaces, as defined in \cite{Ad1975}.
When $s>-\tfrac{1}{2}$, the latter is known to be a Banach algebra. As shown below,
the property extends to $H^{s}_r(\mathbb{R})$, with $s>\tfrac{1}{2}$ and $r\ge0$,  
this fact is crucial for the analysis of Section~\ref{sec4}.
\begin{lemma}\label{lm2.2}
Assume $s>\tfrac{1}{2}$ and $r\ge 0$. Then $H^s_r(\mathbb{R})$ is a Banach algebra, i.e.
for any $f,g\in H^s_r(\mathbb{R})$
\begin{equation}\label{eq2.6}
\|fg\|_{H^s_r(\mathbb{R})} \le c\|f\|_{H^s_r(\mathbb{R})}\|g\|_{H^s_r(\mathbb{R})},
\end{equation}
with $c>0$ independent of $f$ and $g$.
\end{lemma}
\begin{proof}
(a) Using the elementary estimate $|\xi_0+\xi_1|^s \le c\bigl(|\xi_0|^s + |\xi_1|^s\bigr)$,\footnote{
Which holds for all $\xi_0,\xi_1\in\mathbb{R}$ and $s>-1$, with an absolute 
constant $c>0$ that depends on $s$ only.} combined with the standard convolution Young inequality, 
for any two Hermitian functions $\hat{f}, \hat{g}\in L^2_s(\mathbb{R}\pm)\cap L^2(\mathbb{R}_\pm)
=  L^2(\mathbb{R}_\pm, (1+|\xi|^{2s})d\xi) =: \bar{L}^2_s(\mathbb{R}_\pm)$, we have
\[
\|\hat{f}\ast \hat{g}\|_{L^2_{s}(\mathbb{R}_\pm)}\le 
c\bigl(\|\hat{f}\|_{L^1(\mathbb{R}_\pm)}\|g\|_{L^2_s(\mathbb{R}_\pm)}+
\|f\|_{L^2_s(\mathbb{R}_\pm)}\|\hat{g}\|_{L^1(\mathbb{R}_\pm)}\bigr).
\]
By our assumption $s>\tfrac{1}{2}$, hence the direct application of H\"older's inequality yields
\[
\|\hat{f}\|_{L^1(\mathbb{R}_\pm)} 
\le \|(1+|\xi|^{2s})^{-\frac{1}{2}}\|_{L^2(\mathbb{R}_+)}
\|\hat{f}\|_{\bar{L}^2_s(\mathbb{R}_\pm)}
\le c \|\hat{f}\|_{\bar{L}^2_s(\mathbb{R}_\pm)}
\]
and we conclude 
\[
\|\hat{f}\ast \hat{g}\|_{L^2_{s}(\mathbb{R}_\pm)}
\le c \|\hat{f}\|_{\bar{L}^2_s(\mathbb{R}_\pm)}
\|\hat{g}\|_{\bar{L}^2_s(\mathbb{R}_\pm)}.
\]

(b) We let 
\[
\bar{L}^{2,r}_s(\mathbb{R}_\pm) 
:= \mathcal{J}_{\mp}^r[\hat{L}^2_s(\mathbb{R}_\pm)]
= L^{2,r}_s(\mathbb{R}_\pm)\cap L^{2,r}_0(\mathbb{R}_\pm)
\]
and observe that  $\bar{L}^{2,r_1}_s(\mathbb{R}_\pm) 
\hookrightarrow \bar{L}^{2,r_0}_s(\mathbb{R}_\pm) $, whenever $0\le r_0\le r_1$ 
(see \cite[formula (21)]{BaPaPoSh2014}). 
By definition, $\mathcal{P}_+ + \mathcal{P}_- = \mathcal{I}$, where $\mathcal{I}$ is the 
identity operator. Therefore,
\[
\mathcal{P}_+[\kappa_{r,\ell}^- fg] = 
\mathcal{P}_+[\kappa_{\frac{r}{2},\ell}^- f]\mathcal{P}_+[\kappa_{\frac{r}{2},\ell}^- g]
+\mathcal{P}_+[\kappa_{\frac{r}{2},\ell}^- f]\mathcal{P}_-[\kappa_{\frac{r}{2},\ell}^- g]
+\mathcal{P}_-[\kappa_{\frac{r}{2},\ell}^- f]\mathcal{P}_+[\kappa_{\frac{r}{2},\ell}^- g].
\]
Finally, 
$\kappa_{\frac{r}{2},\ell}^- = \sum_{i=0}^{\frac{r}{2}} 
\tbinom{r/2}{\ell} (2i\ell)^{\frac{r}{2}-i}\kappa_{i,\ell}^+$,
provided $\frac{r}{2}$ is a positive integer. These facts, combined with part (a) of the proof, 
yield the bound
\begin{align}
\nonumber
\|\hat{f}\ast\hat{g}\|_{L^{2,r}_\alpha(\mathbb{R}_\pm)} 
&\le c \sum_{i,j=0}^{\frac{r}{2}}
\|\hat{f}\|_{\bar{L}^{2,i}_s(\mathbb{R}_\pm)}
\|\hat{g}\|_{\bar{L}^{2,j}_s(\mathbb{R}_\pm)}\\
\label{eq2.7}
&\le c \|\hat{f}\|_{\bar{L}^{2,\frac{r}{2}}_s(\mathbb{R}_\pm)}
\|\hat{g}\|_{\bar{L}^{2,\frac{r}{2}}_s(\mathbb{R}_\pm)},\quad \tfrac{r}{2}\in \mathbb{N}.
\end{align}

(c) We note that  for any $s>-\tfrac{1}{2}$,
$w = 1+|\xi|^{2s}\in A^{\loc}_{+,2}(\mathbb{R}_+)\cap A^{\loc}_\infty(\mathbb{R})$ (see
Appendix~\ref{sec7}). Hence, by Corollary~\ref{lm7.8} in Appendix~\ref{sec7},
\[
[\bar{L}^{2,r_0}_s(\mathbb{R}_+), \bar{L}^{2,r_1}_s(\mathbb{R}_+)]_\theta
=\bar{L}^{2,(1-\theta)r_0+\theta r_1}_s(\mathbb{R}_+),
\]
$\theta\in(0,1)$, $r_0,r_1\ge 0$, $s>-\tfrac{1}{2}$.
Viewing the convolution product in the Fourier space as a bilinear map from 
$\bar{L}^{2,\frac{r}{2}}_s(\mathbb{R}_+)\times \bar{L}^{2,\frac{r}{2}}_s(\mathbb{R}_+)$ 
to $L^{2,r}_s(\mathbb{R}_+)$, $s>\tfrac{1}{2}$, 
$r\ge 2$ and making use of the classical multilinear complex interpolation theorem of A. Calderon 
(see e.g. \cite[Theorem~4.4.1]{BeLo1976}), we infer from \eqref{eq2.7}
\[
\|\hat{f}\ast\hat{g}\|_{L_s^{2,r}(\mathbb{R}_+)} 
\le c \|\hat{f}\|_{\bar{L}^{2,\frac{r}{2}}_s(\mathbb{R}_+)}
\|\hat{g}\|_{\bar{L}^{2,\frac{r}{2}}_s(\mathbb{R}_+)},\quad s >\tfrac{1}{2},\quad r\ge 2.
\]
By virtue of \cite[formula (21)]{BaPaPoSh2014},
\[
\|\hat{f}\|_{\bar{L}^{2,\frac{r}{2}}_s(\mathbb{R}_\pm)}
\le \|\hat{f}\|_{L^{2,\frac{r}{2}}_0(\mathbb{R}_\pm)} + 
\|\hat{f}\|_{L^{2,\frac{r}{2}}_s(\mathbb{R}_\pm)}
\le c \|\hat{f}\|_{L^{2,r}_s(\mathbb{R}_\pm)},\quad 0\le s\le \tfrac{r}{2},
\]
while the direct application of the convolution Young inequality in the Fourier space, followed by 
\cite[formula (21)]{BaPaPoSh2014}, for all $s>-\tfrac{1}{2}$ and $r\ge 0$ gives
\begin{align*}
\|fg\|_{L^2(\mathbb{R})} &\le c \bigl(\|f\|_{L^2(\mathbb{R})}\|\hat{g}\|_{\bar{L}^2_s(\mathbb{R}_\pm)} 
+ \|\hat{f}\|_{\bar{L}^2_s(\mathbb{R}_\pm)} \|g\|_{L^2(\mathbb{R})}\bigr)\\
&\le c\bigl(\|f\|_{L^2(\mathbb{R})}\|g\|_{L^2(\mathbb{R})} 
+ \|\hat{f}\|_{L^{2,r}_s(\mathbb{R}_\pm)}\|\hat{g}\|_{L^{2,r}_s(\mathbb{R}_\pm)}\bigr).
\end{align*}
Combining the last three inequalities, we conclude 
that \eqref{eq2.6} holds, with $\tfrac{1}{2}<s\le \tfrac{r}{2}$ and $r\ge 2$. 

(d) To complete the proof, we remark that in the standard non-weighted Sobolev settings ($r=0$), 
\eqref{eq2.6} holds for any $s>\tfrac{1}{2}$, see \cite{Ad1975}. 
Hence, the interpolation identity \eqref{eq2.3} 
and part (c) of the proof, combined together, yield \eqref{eq2.6} for any $r\ge 0$.
\end{proof}

\section{Continuous and discrete MTC approximations}\label{sec3}

The Malmquist-Take\-na\-ka-Christov functions $\{\phi_n\}_{n\ge0}$ are defined 
as Fourier preimages of the classical Laguerre functions.\footnote{For an alternative definition, 
and historical remarks see 
\cite{IsWe2019, Wid1994, Wid1995} and references therein.}
That is, for $k\ge 0$, we have
\begin{subequations}\label{eq3.1}
\begin{align}
\label{eq3.1a}
&\mathcal{F}[\phi_{2k}](\xi) = \hat{\phi}_{2k}(\xi)
=\tfrac{\sqrt{\ell}}{\sqrt{2}}\varphi^{0,\ell}_k(\xi),\quad k\geq0, \\ 
\label{eq3.1b}
&\mathcal{F}[\phi_{2k+1}](\xi) = \hat{\phi}_{2k+1}(\xi)
=-i\tfrac{\sqrt{\ell}}{\sqrt{2}}\sgn(\xi) \varphi^{0,\ell}_k(\xi), \quad k\geq0,
\end{align}
\end{subequations}
where 
\begin{subequations}\label{eq3.2}
\begin{equation}\label{eq3.2a}
\varphi_k^{s,\ell}(\xi) = e^{-\frac{\ell|\xi|}{2}}L_k^{(s)}(\ell|\xi|),\quad k\geq0,\quad \ell>0
\end{equation} 
and $L_k^{(s)}(\cdot)$ are the standard generalized Laguerre polynomials \cite{AS1964}. 
Note  that for $s>-1$, the collection $\{\varphi_k^{s,\ell}\}_{k\geq0}$ provides a complete orthogonal 
basis in the weighted space $L^2_{\frac{s}{2}}(\mathbb{R}_+)$. In particular,
\begin{equation}\label{eq3.2b}
\langle\varphi_k^{s,\ell},\varphi_m^{s,\ell}\rangle_{L^2_{\frac{s}{2}}(\mathbb{R}_+)} = 
\int_{\mathbb{R}_+}\varphi_k^{s,\ell}(\xi)\varphi_m^{s,\ell}(\xi)\xi^s d\xi 	
= \tfrac{1}{\ell^{s + 1}}\tfrac{\Gamma(n+s+1)}{\Gamma(n+1)}\delta_{km},\quad k,m\ge 0.
\end{equation} 
\end{subequations}

Straightforward calculations show that 
\begin{subequations}\label{eq3.3}
\begin{align}
\label{eq3.3a}
&\phi_{2k}(x) 
= 2\sqrt{\tfrac{\ell}{\pi}}\Ima \tfrac{(2x+i\ell)^k}{(2x-i\ell)^{k+1}}
 = \tfrac{2}{\sqrt{\pi\ell}}\sin\tfrac{(2k+1)\theta}{2} \sin\tfrac{\theta}{2},\\
 \label{eq3.3b}
&\phi_{2k+1}(x) 
= 2\sqrt{\tfrac{\ell}{\pi}} \Rea\tfrac{(2x+i\ell)^k}{(2x-i\ell)^{k+1}}
= \tfrac{2}{\sqrt{\pi\ell}}\cos\tfrac{(2k+1)\theta}{2} \sin\tfrac{\theta}{2},
\end{align}
\end{subequations}
where $x = \tfrac{\ell}{2}\cot\tfrac{\theta}{2}$, $\theta\in(0,2\pi)$ and $\ell>0$. As evident from 
\eqref{eq3.1} and \eqref{eq3.2}, 
the system $\{\phi_{n} \}_{n\ge 0}$ is a complete orthonormal basis in $L^2(\mathbb{R})$ and 
\[
\langle \phi_{k}, \phi_{m} \rangle_{L^2(\mathbb{R})} = \delta_{km}, \quad k,m\ge 0.
\]

In context of spectral methods, functions $\phi_{n}$, $n\ge 0$, were discovered by C.I. Christov \cite{Cri1982} 
in an attempt to obtain a computational basis that behaves well with respect to the product of its members.
In particular, the following holds
\begin{subequations}\label{eq3.4}
\begin{align}
\label{eq3.4a}
&\phi_{2k}\phi_{2m} 
= \tfrac{1}{2\sqrt{\pi\ell}}\bigl(\phi_{2(k+m)}-\phi_{2(k+m)+2} + \phi_{2(m-k)}-\phi_{2(m-k)-2}\bigr),\\
\label{eq3.4b}
&\phi_{2k+1}\phi_{2m+1} 
= \tfrac{1}{2\sqrt{\pi\ell}}\bigl(-\phi_{2(k+m)}+\phi_{2(k+m)+2} + \phi_{2(m-k)}-\phi_{2(m-k)-2}\bigr),\\
\label{eq3.4c}
&\phi_{2k}\phi_{2m+1} 
= \tfrac{1}{2\sqrt{\pi\ell}}\bigl(\phi_{2(k+m)+1}-\phi_{2(k+m)+3} + \phi_{2(m-k)+1}-\phi_{2(m-k)-1}\bigr).
\end{align}
\end{subequations}
The system $\{\phi_{n} \}_{n\ge 0}$ has a number of attractive computational features, e.g. in view of \eqref{eq3.3}, 
the MTC functions are connected with the trigonometric basis and hence direct and inverse spectral 
transforms can be computed efficiently via Fast Fourier Transform (FFT) algorithm 
\cite{Boyd1990, Cri1982, IsWe2019, Wid1994, Wid1995}. Differentiation and 
computing of the Hilbert transform are also easy \cite{IsWe2019, Wid1994, Wid1995}
\begin{subequations}\label{eq3.5}
\begin{align}
\label{eq3.5a}
&\tfrac{d}{dx}\phi_{2k}  = \tfrac{k+1}{\ell}\phi_{2k+3} 
- \tfrac{2k+1}{\ell}\phi_{2k+1} + \tfrac{k}{\ell}\phi_{2k-1},\\
\label{eq3.5b}
&\tfrac{d}{dx}\phi_{2k+1}  = -\tfrac{k+1}{\ell}\phi_{2k+2} 
+ \tfrac{2k+1}{\ell}\phi_{2k} - \tfrac{k}{\ell}\phi_{2k-2},\\
\label{eq3.5c}
&\mathcal{H}[\phi_{2k}] = \phi_{2k+1},\quad \mathcal{H}[\phi_{2k+1}] = -\phi_{2k}.
\end{align}
\end{subequations}
In the context of the Benjamin equation, identity \eqref{eq3.5c} is particularly important.\footnote{
For the closely related Benjamin-Ono equation, this property is used explicitly in \cite{BoXu2011, BoXu2012}.}

As far as we are aware, the only rigorous approximation result related to the MTC basis is 
the geometric convergence rate of the continuous MTC-Fourier series for functions analytic in the exterior 
of a neighborhood of $\{i,-i\}$ in $\mathbb{C}$ (see \cite{Boyd1990, Wid1995}, the discussion in 
\cite{IsWe2019} and references therein). Unfortunately, in context of differential equations (and in particular 
of \eqref{eq1.1}) the result is not very informative. In the sequel, we derive several alternative error bounds 
directly in $H^s_r(\mathbb{R})$ settings. The estimates form a necessary theoretical background for the 
convergence analysis of an MTC pseudo-spectral scheme, presented in Section~\ref{sec4}. 

\subsection{Projection errors}
Let $n$ be a positive integer, $\mathbb{P}_n$ be the finite dimensional linear space spanned by 
$\{\phi_k(x)\}_{k=0}^n$, $x\in\mathbb{R}$ and $\hat{\mathbb{P}}_n$ be the finite 
dimensional space spanned by $\bigl\{e^{-\tfrac{\ell|\xi|}{2}}x^\xi\bigr\}_{k=0}^n$, $\xi\ge 0$.
In connection,  with $\mathbb{P}_n$ and $\hat{\mathbb{P}}_n$, we define two families of orthogonal projectors
$\mathcal{P}_n:L^2(\mathbb{R})\to \mathbb{P}_n$ and 
$\hat{\mathcal{P}}_n^{s}:L^2_{\frac{s}{2}}(\mathbb{R}_+)\to \hat{\mathbb{P}}_n$, $s>-1$, $n>0$:
\begin{align*}
&\mathcal{P}_n[f] = \sum_{k=0}^n \phi_k \hat{f}_k,\quad \hat{f}_k = \langle f, \phi_k \rangle_{L^2(\mathbb{R})},\\
&\hat{\mathcal{P}}_n^s[f] = \sum_{k=0}^n \tfrac{\ell^{s+1}\Gamma(k+1)}{\Gamma(k+s+1)}
\varphi_k^{s,\ell} \hat{f}_k^{s,\ell},\quad 
\hat{f}_k^{s,\ell} = \langle f, \varphi_k^{s,\ell} \rangle_{L^2_{\frac{s}{2}}(\mathbb{R}_+)}.
\end{align*}
By virtue of \eqref{eq2.4} and \eqref{eq3.1}, for real valued functions we have
\begin{equation}\label{eq3.6}
\|(\mathcal{I} - \mathcal{P}_n)[f]\|_{H^s_r(\mathbb{R})}^2 = 
\bigl\|(\mathcal{I} - \hat{\mathcal{P}}^0_{\lceil \frac{n}{2}\rceil})[\hat{f}]\bigr\|_{L^2(\mathbb{R}_+)}^2
+\bigl\|(\mathcal{I} - \hat{\mathcal{P}}^0_{\lceil \frac{n}{2}\rceil})
[\hat{f}]\bigr\|_{L^{2,r}_s(\mathbb{R}_+)}^2.
\end{equation}
A comprehensive discussion of the Laguerre-type projectors $\hat{\mathcal{P}}_{n}^s$, 
$s>-1$, is found in \cite{BaPaPoSh2014}. In particular, for $s_0,s_1>-\tfrac{1}{2}$ and $r_0,r_1\ge 0$,
Theorems 1 and 2 of \cite{BaPaPoSh2014} give the bounds\footnote{In fact, only the case of $\ell=1$ 
is treated in \cite{BaPaPoSh2014}. Nevertheless, trivial modifications of arguments yield 
\eqref{eq3.7}, \eqref{eq3.8} for any $\ell>0$.}
\begin{subequations}\label{eq3.7}
\begin{align}\label{eq3.7a}
&\bigl\|(\mathcal{I} - \hat{\mathcal{P}}^{s_0}_n)[\hat{f}]\|_{L^{2,r_0}_{\frac{s_1}{2}}(\mathbb{R}_+)}
\le c (\ell n)^{r_0+\frac{s_0-s_1-r_1}{2}} \|\hat{f}\|_{L^{2,r_1}_{\frac{s_0+r_1}{2}}
(\mathbb{R}_+)},\\
\label{eq3.7b}
&s_1 \le s_0+r_0,\quad r_1 \ge s_0-s_1+2r_0
\end{align}
\end{subequations}
and
\begin{subequations}\label{eq3.8}
\begin{align}\label{eq3.8a}
&\bigl\|(\mathcal{I} - \hat{\mathcal{P}}^{s_0}_n)[\hat{f}]\|_{L^{2,r_0}_{\frac{s_1}{2}}(\mathbb{R}_+)}
\le c (\ell n)^{\frac{s_1-s_0-r_1}{2}} \|\hat{f}\|_{L^{2,r_1}_{\frac{s_0+r_1}{2}}(\mathbb{R}_+)},\\
\label{eq3.8b}
&s_1 \ge s_0+r_0,\quad r_1 \ge s_1-s_0.
\end{align}
\end{subequations}
Combining \eqref{eq3.6}, \eqref{eq3.7}, \eqref{eq3.8} and \eqref{eq2.1}, we have 
\begin{lemma}\label{lm3.1}
Assume $s>-1$ and $r_0,r_1\ge 0$. Then
\begin{subequations}\label{eq3.9}
\begin{align}
\label{eq3.9a}
&\|(\mathcal{I} - \mathcal{P}_n)[f]\|_{H^s_{r_0}(\mathbb{R})}  
\le c \bigl(\tfrac{\ell n}{2}\bigr)^{r_0-r_1-s}\|f\|_{H^{r_1}_{2r_1}(\mathbb{R})},
\quad -\tfrac{1}{2}<s\le \tfrac{r_0}{2}\le \tfrac{s+r_1}{2},\\
\label{eq3.9b}
&\|(\mathcal{I} - \mathcal{P}_n)[f]\|_{H^s_{r_0}(\mathbb{R})}  
\le c \bigl(\tfrac{\ell n}{2}\bigr)^{s-r_1}\|f\|_{H^{r_1}_{2r_1}(\mathbb{R})},
\quad 0\le \tfrac{r_0}{2}\le s\le r_1.
\end{align}
\end{subequations}
with a constant $c>0$ independent of $n$ and/or $f$.
\end{lemma}

Lemma~\ref{lm3.1} provides a complete description of the MTC projection errors 
in  $H^s_r(\mathbb{R})$ settings.  In particular, it explains a peculiar disparity in the 
asymptotic of the MTC-Fourier coefficients of closely related holomorphic functions 
$f(x)$, $g(x) = e^{i\xi_0 x} f(x)$, $\xi_0\in\mathbb{R}$, see e.g. examples and discussion 
in \cite{IsWe2019, Wid1994, Wid1995}. 

By virtue of Lemma~\ref{lm3.1}, $|\hat{f}_k|\to 0$ spectrally (faster than any inverse power of $k$), 
provided $\hat{f}(\xi)$ is smooth in $\mathbb{R}_\pm$ and decreases faster than any inverse power 
of $|\xi|$ at infinity.  Since $\hat{g}(\xi) = \hat{f}(\xi-\xi_0)$, the latter condition is violated if 
$\hat{f}(\xi)$ has an integrable singularity at the origin. This is particularly the case when
$f(x)$ is rational, with poles in the upper and lower complex half planes. 

\subsection{Interpolation errors}
Operators $\mathcal{P}_n$ are hard to use in practice
as the integrals of the form $\langle f, \phi_n\rangle_{L^2(\mathbb{R})}$ are impossible to compute 
in most realistic applications. The practical approach consists in replacing 
the inner products with quadratures. In the no boundaries setting of the real line $\mathbb{R}$, 
it is natural to use Gaussian quadratures. The quadrature approximation leads to a rational interpolation 
process, whose properties are briefly discussed below.

For $n=2p-1$, we let
\begin{subequations}\label{eq3.10}
\begin{align}
\label{eq3.10a}
&\langle f,\phi_k\rangle\approx 
\bar{f}_k=\tfrac{\pi}{4\ell p}\sum_{m=0}^{2p-1}(\ell^2+4x_m^2)\phi_k(x_m)f(x_m), \\
\label{eq3.10b}
&x_m=\tfrac{\ell}{2}\cot\bigl(\tfrac{2m+1}{4p}\pi\bigr), \quad 0\le m \le 2p-1.
\end{align}
\end{subequations}
The discrete inner product \eqref{eq3.10} is exact, provided $f \in \mathbb{P}_{n}$.
In practice, we use the discrete spectral coefficients $\bar{f}_k$ and approximate $f$ by 
\begin{equation}\label{eq3.11}
\mathcal{I}_n[f]=\sum_{k=0}^n\bar{f}_k\phi_{k}.
\end{equation}
Directly from \eqref{eq3.3}, \eqref{eq3.10} and \eqref{eq3.11}, it follows that 
\begin{equation}\label{eq3.12}
\mathcal{I}_n[f](x_m)=f(x_m), \quad 0\le m\le 2p-1,
\end{equation}
i.e. $\mathcal{I}_n[\cdot]$ is an interpolation operator. 

Computational properties of $\mathcal{I}_n$ are very similar 
to those of rational Gauss-Chebyshev interpolants, discussed in \cite{SPG2017} and the
generalized Gauss-Laguerre interpolants of \cite{BaPaPoSh2014}. In particular, 
we have
\begin{lemma}\label{lm3.2}
Assume $f \in \mathbb{P}_n$,  $s > -\tfrac{1}{2}$ and $r\ge0$. Then
\begin{equation}\label{eq3.13}
\|f\|_{H^{s}_r(\mathbb{R})} \leq c \bigl(\tfrac{n}{2\ell}\bigr)^{r+|s|-\min\{0,2s\}} 
\|f\|_{L^2(\mathbb{R})},
\end{equation}	 
with a constant $c>0$ independent of $n$ and/or $f$.
\end{lemma}
\begin{proof}
Since $f\in\mathbb{P}_n$, $n=2p-1$, we have $\hat{f}\in \hat{\mathbb{P}}_p$, $\xi\in\mathbb{R}_+$.
In \cite[Lemma 6]{BaPaPoSh2014}, it is shown that for such functions
\begin{align*}
&\|\hat{f}\|_{L^{2,r}_{s}(\mathbb{R}_+)} 
\le c (\ell p)^{r-\min\{0,2s\}}\|\hat{f}\|_{L^2_{s}(\mathbb{R}_+)}, \\
&\|\hat{f}\|_{L^{2}_{s}(\mathbb{R}_+)} 
\le c (\ell p)^{|s|}\|\hat{f}\|_{L^2(\mathbb{R}_+)}.
\end{align*} 
In view of \eqref{eq2.4}, these inequalities imply \eqref{eq3.13}.
\end{proof}

\begin{lemma}\label{lm3.3} 
Assume $s>\frac{1}{2}$. Then 
\begin{equation}\label{eq3.14}
\|\mathcal{I}_n\|_{H^{s}(\mathbb{R})\to L^2(\mathbb{R})} \le c
\bigl(\tfrac{\ell n}{2}\bigr),
\end{equation}
with $c>0$ independent of $n$.
\end{lemma}
\begin{proof}
Since the discrete inner product \eqref{eq3.10} is exact for $f \in \mathbb{P}_{n}$, we have
\[
\| \mathcal{I}_n[f]\|^2_{L^2(\mathbb{R})}=\tfrac{\pi}{4\ell p}
\sum_{m=0}^{2p-1}(\ell^2+4x_m^2){f^2}(x_m).
\]
In view of the classical Sobolev embedding \cite{Ad1975}, 
$\| f\|_{L^\infty(\mathbb{R})}\leq c_s \|f\|_{H^{s}(\mathbb{R})}$, $s>\frac{1}{2}$. 
Consequently,
\begin{align*}
\| \mathcal{I}_n[f]\|^2_{L^2(\mathbb{R})}&\le \tfrac{c_s\pi}{4\ell p}
\Bigl[\sum_{m=0}^{2p-1}(\ell^2+4x_m^2)\Bigr]\|f\|^2_{H^{s}(\mathbb{R})}\\
&= c \Bigl[\sum_{m=0}^{2p-1}\tfrac{\ell^2}{\sin^2\bigl(\frac{(2m+1)\pi}{4p}\bigr)}\Bigr]
\|f\|^2_{H^{s}(\mathbb{R})}= c S^2_n\|f\|^2_{H^{s}(\mathbb{R})}.
\end{align*}
In \cite[Lemma 4]{SPG2017} it is shown $S^2_n = 2(2\ell p)^2$. Hence, \eqref{eq3.14} is settled. 
\end{proof}

The interpolation error bounds are obtained combining Lemmas~\ref{lm3.1}-\ref{lm3.3}.
\begin{corollary}\label{lm3.4}
Let $s > -\tfrac{1}{2}$, $r_0\ge0$, $\varepsilon > 0$ and $r_1\ge r_0+|s|$. 
Then,
\begin{equation}\label{eq3.15}
\|(\mathcal{I}-\mathcal{I}_n)[f]\|_{H^{s}_{r_0}(\mathbb{R})} 
\le c \bigl(\tfrac{\ell n}{2}\bigr)^{\frac{3}{2}+\varepsilon+r_0+|s|-\max\{0,2s\}-r_1}
\|f\|_{H^{r_1}_{2r_1}(\mathbb{R})},
\end{equation}
with a constant $c>0$ independent of $n$ and/or $f$.
\end{corollary}
\begin{proof} 
In view of Lemmas~\ref{lm3.2} and \ref{lm3.3}, we have
\begin{align*}
&\|(\mathcal{I}-\mathcal{I}_n)[f]\|_{H^{s}_{r_0}(\mathbb{R})}
\le \|(\mathcal{I}-\mathcal{P}_n)[f]\|_{H^{s}_{r_0}(\mathbb{R})}
+\|\mathcal{I}_n(\mathcal{I}-\mathcal{P}_n)[f]\|_{H^{s}_{r_0}(\mathbb{R})}\\
&\qquad\le \|(\mathcal{I}-\mathcal{P}_n)[f]\|_{H^{s}_{r_0}(\mathbb{R})}
+ c\bigl(\tfrac{\ell n}{2}\bigr)^{r_0+|s|-\min\{0,2s\}}
\|\mathcal{I}_n(\mathcal{I}-\mathcal{P}_n)[f]\|_{L^2(\mathbb{R})}\\
&\qquad\le \|(\mathcal{I}-\mathcal{P}_n)[f]\|_{H^{s}_{r_0}(\mathbb{R})}
+ c\bigl(\tfrac{\ell n}{2}\bigr)^{1+r_0+|s|-\min\{0,2s\}}
\|(\mathcal{I}-\mathcal{P}_n)[f]\|_{H^{\frac{1}{2}+\varepsilon}(\mathbb{R})}.
\end{align*}
Hence, \eqref{eq3.15} is the direct consequence of Lemma~\ref{lm3.1}.
\end{proof}

\section{An MTC collocation scheme}\label{sec4}

To obtain a spatial semi-discretization, for a given $n=2p-1$, $p\in\mathbb{N}$, we approximate 
the automorphism $\mathcal{J}$ by the finite dimensional skew symmetric map 
$\mathcal{J}_{n} = -\mathcal{P}_n \partial_x \mathcal{P}_n:\mathbb{P}_n \to \mathbb{P}_n$ and replace
\eqref{eq1.1} with 
\begin{subequations}\label{eq4.1}
\begin{align}
\label{eq4.1a}
&\bar{u}_t = \mathcal{J}_n \nabla_{\bar{u}} \mathcal{G}_{n}(\bar{u}),\quad \bar{u}(0)=\mathcal{I}_n[u_0],\\
\label{eq4.1b}
&\mathcal{G}_{n}(\bar{u})=\tfrac{1}{2}\int_{\mathbb{R}}
\Bigl(\alpha|\bar{u}|^2 - \beta\bar{u}\mathcal{H}[\mathcal{J}_n\bar{u}]+\gamma|\mathcal{J}_n\bar{u}|^2  
+ \tfrac{2\delta}{3}\bar{u}\mathcal{I}_n[\bar{u}^2]  \Bigr) dx,
\end{align}
\end{subequations}
where $\bar{u}\in \mathbb{P}_n$. Note that if $n=2p-1$, the operator $\mathcal{J}_n$ is non-degenerate. 
This follows from identities \eqref{eq3.5a}-\eqref{eq3.5b} and the fact that the eigenvalues of the differentiation 
matrix $-\mathcal{J}_n$ are given explicitly by $\pm i \tfrac{\xi_n}{\ell}$, $1\le k\le p$, 
where $\xi_k$ are roots of the classical Laguerre polynomial $L_p(x)$ (see the proof of Lemma~\ref{lm4.1} below and 
\cite{Wid1994}). As a consequence, the finite dimensional semi-discrete system \eqref{eq4.1} of ODEs 
is again Hamiltonian. 

By construction, the semi-discrete vector field $\nabla \mathcal{G}_n(\bar{u})$ is smooth 
and hence the initial value problem \eqref{eq4.1a} is locally well-posed. Unfortunately, the only 
conserved quantity\footnote{This is in contrast with the exact classical solutions, where,  in addition to the 
Hamiltonian, the $L^2(\mathbb{R})$ norm is preserved.} $\mathcal{G}_n(\bar{u})$ is indefinite. 
As a consequence, we have insufficient amount of a priori information to establish uniform global bounds 
on the growth rate of the numerical solution $\bar{u}$. 
To alleviate the problem, we proceed indirectly. Instead of estimating $\bar{u}$, we compare it 
to the reference solution $\tilde{u} = \mathcal{P}_n[u]\in\mathbb{P}_n$, where 
$u$ (the exact classical solution to \eqref{eq1.1}) is assumed to be globally defined and regular.\footnote{
The approach is a manifestation of an elementary observation that in the Cauchy problem 
$y' = y^s$, $y(0) = y_0$, $s>1$, the blow up time is inverse proportional to the size of the input data. 
The idea is widely used in numerical analysis, see e.g. \cite{MaQu1988} for an application in the context 
of spectral methods.}

\subsection{Auxiliary estimates}
In our analysis, we make use of three technical estimates. The first one is a discrete analogue 
of the classical Gagliardo-Nirenberg inequality, the second is used to estimate discrete 
power nonlinearities and the last one is an extension of the classical Gronwall's Lemma. 
\begin{lemma}\label{lm4.1}
Let $u\in\mathbb{P}_n$, $n=2p-1$. Then
\begin{subequations}\label{eq4.2}
\begin{align}
\label{eq4.2a}
&\|u_x\|_{L^2(\mathbb{R})} \le c \bigl(\tfrac{\ell n}{2}\bigr)^{\frac{1}{2}}
\|\mathcal{J}_n u\|_{L^2(\mathbb{R})},\\
\label{eq4.2b}
&\|u\|_{L^\infty(\mathbb{R})} \le c \bigl(\tfrac{\ell n}{2}\bigr)^{\frac{1}{4}}
\|u\|_{L^2(\mathbb{R})}^{\frac{1}{2}}\|\mathcal{J}_n u\|_{L^2(\mathbb{R})}^{\frac{1}{2}},
\end{align}
\end{subequations}
where $c>0$ is an absolute constant.
\end{lemma}
\begin{proof}
Identities \eqref{eq3.5} imply
\[
\|u_x\|_{L^2(\mathbb{R})}^2 = \|\mathcal{J}_n u\|^2_{L^2(\mathbb{R})} + 
\tfrac{p^2}{\ell^2}\bigl[|\hat{u}_{2(p-1)}|^2 + |\hat{u}_{2p-1}|^2\bigr],\quad 
\]
where $\hat{u}_k = \langle u, \phi_k \rangle_{L^2(\mathbb{R})}$, $0\le k\le n$.
Our main task is to bound the sum $|\hat{u}_{2(p-1)}|^2+|\hat{u}_{2p-1}|^2$.

Let $\hat{u}_e = (\hat{u}_0,\ldots,\hat{u}_{2(p-1)})^T$ and
$\hat{u}_o = (\hat{u}_1,\ldots,\hat{u}_{2p-1})^T$ be $\mathbb{R}^p$ vectors that 
contain the even and the odd MTC-Fourier coefficients of $u\in\mathbb{P}_n$.
Then, by virtue of \eqref{eq3.5}, the even and the odd MTC-Fourier coefficients of $-\mathcal{J}_n u$ 
are given by $\tfrac{1}{\ell}D \hat{u}_o$ and $-\tfrac{1}{\ell}D\hat{u}_e$, respectively, 
where $D=(d_{ij})\in\mathbb{R}^{p\times p}$ is the symmetric three-diagonal matrix, whose entries 
are given by $d_{ii} = -2i-1$, $d_{i,i+1} = d_{i+1,i} = i$, $0\le i\le p$.  
Using the three-term recurrence formula for the classical Laguerre polynomials $L_n(x)$ (see \cite{AS1964, Wid1994}), 
we find that 
\[
D = Q\Lambda Q^T, \quad \Lambda = \diag(\xi_0,\ldots,\xi_{p-1}),
\quad Q_{ij} = \tfrac{\sqrt{\xi_j}L_i(\xi_j)}{p|L_{p-1}(\xi_j)|},
\]
where $\xi_i$, $0\le i\le p-1$, are the (strictly positive) roots of $L_p(x)$ and that 
matrix $Q\in\mathbb{R}^p$ is orthogonal. 

Let $e_p$ be the standard unit vector in $\mathbb{R}^p$ and $|\cdot|$, $\cdot$ denote 
the usual Euclidean norm and the inner product in $\mathbb{R}^p$. With this notation, we obtain
\begin{align*}
|\hat{u}_{2(p-1)}|^2+|\hat{u}_{2p-1}|^2 &= |e_p\cdot \hat{u}_e|^2+|e_p\cdot \hat{u}_o|^2\\
&\le \ell^2 |\Lambda^{-1} Q^T e_p|^2 \|\mathcal{J}_nu\|_{L^2(\mathbb{R})}^2.
\end{align*}
Note that 
\[
c\xi_i \le \tfrac{(i+1)^2}{p} \le C\xi_i,\quad 0\le i\le p-1,
\]
for some absolute constants $c,C>0$ (see e.g. \cite[formula (2.3.50), p. 141]{Mastr}). Hence, 
\[
|\Lambda^{-1} Q^T e_p|^2 = \tfrac{1}{p^2}\sum_{i=0}^{p-1}\tfrac{1}{\xi_i} \le \tfrac{c}{p}
\]
and \eqref{eq4.2a} follows. Bound \eqref{eq4.2b} follows from \eqref{eq4.2a} and the standard
Gagliardo-Nirenberg inequality.
\end{proof}

\begin{lemma}\label{lm4.2}
Assume $v\in\mathbb{P}_n$, $m>0$, $2\le k\le 5$ and $1\le r\le 2$. Then
\begin{subequations}\label{eq4.3}
\begin{align}
\label{eq4.3a}
\bigl|\langle \mathcal{I}_n[\tilde{u}^m], \mathcal{I}_n[v^k]\rangle_{L^2(\mathbb{R})}\bigr|
& \le c\bigl(\tfrac{\ell n}{2\varepsilon}\bigr)^{\frac{k-2}{6-k}}
\|v\|_{L^2(\mathbb{R})}^{\frac{2(k+2)}{6-k}} 
+\varepsilon (k-2)\|\mathcal{J}_n v\|_{L^2(\mathbb{R})}^2,\\
\label{eq4.3b}
\bigl|\langle \mathcal{I}_n[\tilde{u}^m], \mathcal{I}_n[v^r\mathcal{J}_n v]\rangle_{L^2(\mathbb{R})}\bigr|
&\le c\varepsilon^{-\frac{2}{3-r}}\bigl(\tfrac{\ell n}{2\varepsilon}\bigr)^{\frac{r-1}{3-r}}
\|v\|_{L^2(\mathbb{R})}^{\frac{2(r+1)}{3-r}}
+\varepsilon\|\mathcal{J}_n v\|_{L^2(\mathbb{R})}^2,
\end{align}
\end{subequations}
where $\varepsilon>0$ is arbitrary and $c>0$ depends on $k$ and $\|\tilde{u}\|_{L^\infty(\mathbb{R})}$ only.
\end{lemma}
\begin{proof}
Let $w_i = \tfrac{\pi}{4\ell p}(\ell^2+4x_i^2)$, $0\le i\le 2p-1$, where $x_i$ is defined in \eqref{eq3.10b}.
Since quadrature \eqref{eq3.10a} is exact in $\mathbb{P}_n$ and in view of Lemma~\ref{lm4.1}, we have 
\begin{align*}
\bigl|\langle \mathcal{I}_n[\tilde{u}^m], \mathcal{I}_n[v^{k}]\rangle_{L^2(\mathbb{R})}\bigr|
&\le \sum_{i=0}^{n} w_i |\tilde{u}(x_i)|^m |v(x_i)|^k\\
&\le \|\tilde{u}\|^m_{L^\infty(\mathbb{R})}\|v\|_{L^\infty(\mathbb{R})}^{k-2} \|v\|^2_{L^2(\mathbb{R})}\\
&\le c\|\tilde{u}\|^m_{L^\infty(\mathbb{R})}\|v_x\|_{L^\infty(\mathbb{R})}^{\frac{k-2}{2}} 
\|v\|^{\frac{k+2}{2}}_{L^2(\mathbb{R})}\\
&\le c \bigl(\tfrac{\ell n}{2}\bigr)^{\frac{k-2}{4}}
\|\tilde{u}\|^m_{L^\infty(\mathbb{R})}
\|\mathcal{J}_n v\|_{L^2(\mathbb{R})}^{\frac{k-2}{2}} 
\|v\|^{\frac{k+2}{2}}_{L^2(\mathbb{R})}.
\end{align*}
Hence, Young's inequality, with exponents $\tfrac{4}{k-2}$ and $\tfrac{4}{6-k}$, yields \eqref{eq4.3a}.
The proof of \eqref{eq4.3b} is identical. 
\end{proof}

\begin{lemma}\label{lm4.3}
Let $u\in C[0,T]$ be non-negative. Assume that  
\begin{subequations}\label{eq4.4}
\begin{equation}\label{eq4.4a}
u(t) \le f(t) + a\int_0^t u(s)ds + b\int_0^t (t-s) u(s)ds,\quad t\in [0,T],
\end{equation}
where $a,b>0$ and $f(t)$ is integrable and non-negative.
Then
\begin{equation}\label{eq4.4b}
u(t) \le f(0) + e^{\frac{a+\sqrt{a^2+4b}}{2}t}\int_0^t \bigl[1+\tfrac{b}{2}(t-s)^2\bigr]f(s)ds,
\quad t\in[0,T].
\end{equation}
\end{subequations}
\end{lemma}
\begin{proof}
Let $U(t) = e^{\lambda t}\int_0^t (t-s) u(s)ds$, where $\lambda = -\frac{a+\sqrt{a^2+4b}}{2}$ is 
the negative root of the quadratic equation $\lambda^2 +a \lambda - b=0$. 
Then \eqref{eq4.4a} is equivalent to
\[
U''(t) \le (a+2\lambda)U'(t) + f(t),\quad U(0) = U'(0) = 0, \quad t\in [0,T].
\]
Since $a+2\lambda<0$, integrating twice, we obtain
\[
e^{-\lambda t}U(t) = \int_0^t (t-s) u(s)ds \le e^{-\lambda t}\int_0^t (t-s) f(s)ds.
\]
Upon substitution into \eqref{eq4.4a}, we have
\[
u(t) \le f(t)+ e^{-\lambda t}b\int_0^t (t-s) f(s)ds + a\int_0^t u(s)ds,
\]
which, combined with the standard Gronwall's inequality, gives \eqref{eq4.4b}.
\end{proof}

\subsection{Stability}
Now, we turn to the study of the numerical error $e = \tilde{u} - \bar{u}$. 
Applying operator $\mathcal{P}_n$ to both sides of \eqref{eq1.1}, subtracting \eqref{eq4.1} 
and passing to the quadrature (as we did in Lemma~\ref{lm4.2}), we infer
\begin{subequations}\label{eq4.5}
\begin{align}
\label{eq4.5a}
e_t  &= \mathcal{J}_n\nabla_{e}\bigl[\mathcal{G}_n(-e) + \mathcal{E}_{n}(e, t)
+ \mathcal{D}_n(e,t)\bigr], \quad e(0) = e_0.\\
\label{eq4.5b}
\mathcal{E}_{n}(e, t) &= 2\delta
\bigl\langle \tilde{u}(t), \mathcal{I}_n[e^2]\bigr\rangle_{L^2(\mathbb{R})},\\
\nonumber
\mathcal{D}_{n}(e,t) &= \bigl \langle e, \alpha(\mathcal{I}-\mathcal{P}_n)\bigl[u](t) 
+ \beta\mathcal{H}\bigl[(\partial_x+\mathcal{J}_n)u\bigr](t),\\
\label{eq4.5c}
&-\gamma(\partial_{xx} - \mathcal{J}_n^2)[u](t) + \delta\bigl(u^2(t) - \mathcal{I}_n\bigl[\tilde{u}^2\bigr](t)\bigr) 
\bigr\rangle_{L^2(\mathbb{R})},
\end{align}
\end{subequations}
where $\nabla_e$ denotes the gradient with respect to variable $e$.
Equation \eqref{eq4.5a} is not Hamiltonian. Nevertheless, differentiating and using the 
skew-symmetry of the discrete automorphism $\mathcal{J}_n$, 
we obtain 
\[
\tfrac{d}{dt}\bigl[ \mathcal{G}_n(-e) +  \mathcal{E}_{n}(e, t) + \mathcal{D}_n(e,t)\bigr]
 = \partial_t[ \mathcal{E}_{n}(e, t) + \mathcal{D}_n(e,t)],
\]
which, after integration in time,  gives
\begin{align}
\nonumber
|\mathcal{G}_n(-e)| 
& \le |\mathcal{G}_n(-e_0)| + |\mathcal{E}_{n}(e_0, 0)| + |\mathcal{D}_n(e_0,0)| \\
\nonumber
&+ |\mathcal{E}_{n}(e, t)| + |\mathcal{D}_n(e,t)|  \\
\label{eq4.6}
&+ \int_0^t \bigl|\partial_t[ \mathcal{E}_{n}(e(s), s) + \mathcal{D}_n(e(s),s)]\bigr|ds.
\end{align}
We use \eqref{eq4.6} to control the $L^2(\mathbb{R})$ norm of $\mathcal{J}_n e$.
\begin{lemma}\label{lm4.4}
Let $\gamma>0$, $0<\varepsilon<\frac{\gamma}{4}$ and
\begin{subequations}\label{eq4.7}
\begin{equation}\label{eq4.7a}
\|e\|_{L^2(\mathbb{R})} \le \bigl(\tfrac{\ell n}{2\varepsilon}\bigr)^{-\frac{1}{4}},
\end{equation}
in some interval $[0,T]$. Then for each $t\in [0,T]$, we have
\begin{align}
\nonumber
\|\mathcal{J}_n e\|^2_{L^2(\mathbb{R})} &\le c\Bigl(
|\mathcal{G}_n(-e_0)| + |\mathcal{E}_{n}(e_0, 0)| + |\mathcal{D}_n(e_0,0)|\\
\nonumber
&+ \|e\|^2_{L^2_n(\mathbb{R})}
+\int_0^t \|e\|^2_{L^2_n(\mathbb{R})}ds\\
\label{eq4.7b}
&+\|\nabla_e \mathcal{D}_n(e,t)\|^2_{L^2(\mathbb{R})}
+\int_0^t \|\nabla_e \partial_t \mathcal{D}_n(e,t)\|^2_{L^2(\mathbb{R})}ds\Bigr),
\end{align}
\end{subequations}
where $c>0$ depends on $\alpha$, $\beta$, $\gamma$, $\delta$, $\|\tilde{u}\|_{L^\infty(\mathbb{R})}$ 
and $\|\tilde{u}_t\|_{L^\infty(\mathbb{R})}$ only.
\end{lemma}
\begin{proof}
We bound each term in \eqref{eq4.6} separately.
First, we use the Cauchy-Schwartz inequality, unitarity of $\mathcal{H}$ and Lemma~\ref{lm4.2} to obtain
\[
\mathcal{G}_n(-e) \ge 
\bigl(\tfrac{\gamma}{2}-\varepsilon\bigr)\|\mathcal{J}_n e\|^2_{L^2(\mathbb{R})}
- c_1\|e\|^2_{L^2(\mathbb{R})}
- c_1\bigl(\tfrac{\ell n}{2\varepsilon}\bigr)^{\frac{1}{3}} \|e\|^{\frac{10}{3}},
\] 
with $c_1>0$ that depends on the parameters $\alpha$, $\beta$, $\gamma$ and $\delta$ only.
Using Lemma~\ref{lm4.2}, we have also
\begin{align*}
&|\mathcal{E}_n(e, t)| \le c_2\|e\|^2_{L^2(\mathbb{R})},\\
&|\partial_t \mathcal{E}_n(e, t)| \le c_3\|e\|^2_{L^2(\mathbb{R})},
\end{align*}
where $c_2>0$ depends on $\delta$ and $\|\tilde{u}\|_{L^\infty(\mathbb{R})}$ and $c_3>0$ depends 
on $\delta$ and $\|\tilde{u}_t\|_{L^\infty(\mathbb{R})}$ only.
The quantity $\mathcal{D}_n(e,t)$ is linear in $e$ and by the Minkowski inequality,
\begin{align*}
&2|\mathcal{D}_n(e,t)| \le \|e\|^2_{L^2(\mathbb{R})}
+\|\nabla_e \mathcal{D}_n(e,t)\|^2_{L^2(\mathbb{R})},\\ 
&2|\partial_t \mathcal{D}_n(e,t)| \le \|e\|_{L^2(\mathbb{R})}^2
+\|\nabla_e \partial_t \mathcal{D}_n(e,t)\|_{L^2(\mathbb{R})}^2.
\end{align*}
Hence \eqref{eq4.7b} is the direct consequence of the above bounds, \eqref{eq4.6} and 
assumption \eqref{eq4.7a}.
\end{proof}

We remark that the assumption $\gamma>0$ appearing in Lemma~\ref{lm4.4} is not restrictive, 
for if $\gamma<0$ one can use $-\mathcal{G}_n(\cdot)$ instead of $\mathcal{G}_n(\cdot)$.

\begin{theorem}[Stability]\label{lm4.5}
Assume that for some fixed $C>0$, $T>0$ and $\epsilon>0$, 
\begin{subequations}\label{eq4.8}
\begin{align}\label{eq4.8a}
&\max\{\|\tilde{u}\|_{L^\infty([0,T]\times \mathbb{R})}, 
\|\tilde{u}_t\|_{L^\infty([0,T]\times \mathbb{R})}\} < C,\\
\nonumber
&\|e_0\|_{L^2(\mathbb{R})} + |\mathcal{G}_n(-e_0)|^{\frac{1}{2}}
+ |\mathcal{E}_{n}(e_0, 0)|^{\frac{1}{2}} + |\mathcal{D}_n(e_0,0)|^{\frac{1}{2}}\\
\nonumber
&\qquad\qquad+ \|\nabla_e \mathcal{D}_n(e,t)\|_{L^2([0,T]\times\mathbb{R})}
+ \|\partial_t \nabla_e\mathcal{D}_n(e,t)\|_{L^2([0,T]\times\mathbb{R})}\\
\label{eq4.8b}
&\qquad\qquad = \mathcal{O}\Bigl(\bigl(\tfrac{\ell n}{2}\bigr)^{-\frac{1+\epsilon}{4}}\Bigr),
\end{align}
uniformly for large values of $n=2p-1>0$. 
Then there exists $c>0$, that depends on 
$C$, $T$ and parameters $\alpha$, $\beta$, $\gamma$ and $\delta$ of \eqref{eq1.1} only, such that 
\begin{align}
\nonumber
\|e\|_{C([0,T], L^2(\mathbb{R}))}
&\le c\Bigl( \|e_0\|_{L^2(\mathbb{R})}
+ |\mathcal{G}_n(-e_0)|^{\frac{1}{2}}
+ |\mathcal{E}_{n}(e_0, 0)|^{\frac{1}{2}} + |\mathcal{D}_n(e_0,0)|^{\frac{1}{2}}\\
\label{eq4.8c}
&+  \|\nabla_e \mathcal{D}_n(e,t)\|_{L^2([0,T]\times\mathbb{R})}
+ \|\nabla_e\partial_t \mathcal{D}_n(e,t)\|_{L^2([0,T]\times\mathbb{R})}
\Bigr),
\end{align}
\end{subequations}
for all sufficiently large values of $n>0$.
\end{theorem}
\begin{proof}
(a) We multiply both sides of \eqref{eq4.4a} by $e$, integrate with respect to $x$ over $\mathbb{R}$ 
and take into account the skew-symmetry of the automorphism $\mathcal{J}_n$. This gives
\begin{align*}
\tfrac{1}{2}\tfrac{d}{dt}\|e\|^2_{L^2(\mathbb{R})} &= 
- \langle \mathcal{J}_n[e], \mathcal{I}_n[e^2]\rangle_{L^2(\mathbb{R})}\\
&- \langle \mathcal{J}_n[e], \nabla_e \mathcal{E}_n(e,t)\rangle_{L^2(\mathbb{R})}
- \langle \mathcal{J}_n[e], \nabla_e \mathcal{D}_n(e,t)\rangle_{L^2(\mathbb{R})}.
\end{align*}
Lemma~\ref{lm4.2} and the Cauchy-Schwartz inequality give the bounds
\begin{align*}
|\langle \mathcal{J}_n[e], \mathcal{I}_n[e^2]\rangle_{L^2(\mathbb{R})}|
&\le \varepsilon\|\mathcal{J}_n e\|_{L^2(\mathbb{R})}^2
+c_1\bigl(\tfrac{\ell n}{2\varepsilon}\bigr)\varepsilon^{-2}\|e\|_{L^2(\mathbb{R})}^{6},\\ 
|\langle \mathcal{J}_n[e], \nabla_e \mathcal{E}_n(e,t)\rangle_{L^2(\mathbb{R})}| 
&\le \|\mathcal{J}_n e\|^2_{L^2(\mathbb{R})}+ c_2\|e\|^2_{L^2(\mathbb{R})},\\
|\langle \mathcal{J}_n[e], \nabla_e \mathcal{D}_n(e,t)\rangle_{L^2(\mathbb{R})}|
&\le \|\mathcal{J}_n e\|^2_{L^2(\mathbb{R})} + 
c_3\|\nabla_e \mathcal{D}_n(e,t)\|^2_{L^2(\mathbb{R})},
\end{align*}
where $c_1,c_3>0$ are absolute constants and $c_2>0$ depends on 
$\|\tilde{u}\|_{L^\infty([0,T]\times \mathbb{R})}$ only.

(b) In view of \eqref{eq4.8b} and local continuity of $\|e\|^2_{L^2(\mathbb{R})}$, we see that 
\eqref{eq4.7a} holds locally in some nonempty closed interval $[0,\tau_0]$, $0<\tau_0\le T$.
Therefore, combining our estimates from part (a) of the proof and using \eqref{eq4.7a} with 
$\varepsilon = \mathcal{O}(1)$, $0<\varepsilon<\tfrac{\gamma}{2}$, we conclude that 
the following holds
\begin{align*}
\tfrac{1}{2}\tfrac{d}{dt}\|e\|^2_{L^2(\mathbb{R})} 
&\le \bigl(c_1+c_2\bigr)\|e\|^2_{L^2(\mathbb{R})}\\
&+ (2+\varepsilon) \|\mathcal{J}_n e\|_{L^2(\mathbb{R})}
+c_3\|\nabla_e \mathcal{D}_n(e,t)\|^2_{L^2(\mathbb{R})},
\end{align*} 
uniformly in $[0,\tau_0]$. Integrating the last formula with respect to time and combining 
the result with Lemma~\ref{lm4.4}, we obtain
\begin{align}
\nonumber
\|e\|^2_{L^2(\mathbb{R})}
&\le \|e_0\|^2_{L^2(\mathbb{R})}
+c_4t\bigl(|\mathcal{G}_n(-e_0)| + |\mathcal{E}_{n}(e_0, 0)| + |\mathcal{D}_n(e_0,0)|\bigr)\\
\nonumber
&+c_4\int_0^t (1+t-s)\|e\|^2_{L^2(\mathbb{R})}ds
+c_4\int_0^t \|\nabla_e \mathcal{D}_n(e,t)\|^2_{L^2(\mathbb{R})}ds\\
\label{eq4.9}
&+c_4\int_0^t (t-s) \|\nabla_e\partial_t \mathcal{D}_n(e,t)\|^2_{L^2(\mathbb{R})}ds,
\end{align} 
where $t\in[0,\tau_0]$ and $c>0$ depends on $C>0$ and parameters $\alpha$, $\beta$, $\gamma$ 
and $\delta$ of the model \eqref{eq1.1} only. 

(c) Inequality \eqref{eq4.9} falls in the scope of Lemma~\ref{lm4.3}, 
hence, definitely \eqref{eq4.8c} holds in the small interval $[0,\tau_0]$. Furthermore, 
from the same Lemma~\ref{lm4.3}, it follows that the constant $c>0$ in \eqref{eq4.8c} 
behaves like $c'(1+\tau_0^{\frac{3}{2}})e^{c'' \tau_0}$, where $c',c''>0$ are independent 
of $n>0$ and $\tau_0$. The observation implies that $\|e(\tau_0)\|^2_{L^2(\mathbb{R})} 
= \mathcal{O}\Bigl(\bigl(\tfrac{\ell n}{2}\bigr)^{-\frac{1+\epsilon}{4}}\Bigr)$,
i.e for $n>0$ sufficiently large, \eqref{eq4.7a} is satisfied at the endpoint $\tau_0$. 
In view of the last fact and by continuity of $\|e\|^2_{L^2(\mathbb{R})}$, we conclude that 
\eqref{eq4.9} can be extended to a larger interval $[0,\tau_1]$,  $0<\tau_0<\tau_1\le T$, 
without increasing the size of the constant $c_4>0$. 

The assertion of Theorem~\ref{lm4.5} follows from the standard continuation argument. 
Repeating the continuation step described above inductively, we construct an ascending sequence 
$0<\tau_0<\tau_1<\cdots \le T$, such that \eqref{eq4.9} (with $c_4>0$ being fixed) holds in each 
$[0,\tau_i]$, $i\ge 0$. Assuming $\tau^\ast = \sup \tau_i <T$, we arrive at the contradiction; 
for if $\tau^\ast<T$, the continuation step extends \eqref{eq4.9} beyond the interval $[0,\tau^\ast]$.
\end{proof}

\subsection{Consistency and convergence}

In what follows, we use the results of Sections~\ref{sec2} and \ref{sec3} to demonstrate that 
assumptions \eqref{eq4.8a}, \eqref{eq4.8b} are satisfied, provided the exact solution $u$ is sufficiently regular.
We begin with \eqref{eq4.8a}.
\begin{lemma}\label{lm4.6}
Assume $u, u_t\in L^\infty([0,T], H^s_{2s}(\mathbb{R}))$, $s>1$. Then 
\begin{subequations}\label{eq4.10}
\begin{align}
\label{eq4.10a}
&\|\tilde{u}\|_{L^\infty([0,T]\times \mathbb{R})} \le c\bigl[1+\bigl(\tfrac{\ell n}{2}\bigr)^{1-s}\bigr]
\|u\|_{L^\infty([0,T], H^s_{2s}(\mathbb{R}))},\\
\label{eq4.10b}
&\|\tilde{u}_t\|_{L^\infty([0,T]\times \mathbb{R})} \le c\bigl[1+\bigl(\tfrac{\ell n}{2}\bigr)^{1-s}\bigr]
\|u_t\|_{L^\infty([0,T], H^s_{2s}(\mathbb{R}))},
\end{align}
\end{subequations}
where $c>0$ is an absolute constant.
\end{lemma}
\begin{proof}
By the standard Gagliardo-Nirenberg inequality, 
\[
\|\tilde{u}\|_{L^\infty(\mathbb{R})}^2 \le c \|\tilde{u}\|_{L^2(\mathbb{R})} \|\tilde{u}_x\|_{L^2(\mathbb{R})}.
\]
Since the Christov functions form a complete orthogonal basis in $L^2(\mathbb{R})$, we have 
$\|\tilde{u}\|_{L^2(\mathbb{R})} = \|\mathcal{P}_n[u]\|_{L^2(\mathbb{R})}\le \|u\|_{L^2(\mathbb{R})}$.
To bound the norm of $\tilde{u}_x$, we write
\[
\|\tilde{u}_x\|_{L^2(\mathbb{R})} \le \|u_x\|_{L^2(\mathbb{R})} + 
\|(\mathcal{I}-\mathcal{P}_n)[u]\|_{H^1(\mathbb{R})},
\]
and apply Lemma~\ref{lm3.1}. This gives \eqref{eq4.10a}. Bound \eqref{eq4.10b} follows along the same lines.
\end{proof}

Next, we show that each term in \eqref{eq4.8b} is small.
\begin{lemma}\label{lm4.7}
Assume $u, u_t\in L^\infty([0,T], H^s_{2s}(\mathbb{R}))$ and $s\ge 2$ and 
$\epsilon>0$. 
Then
\begin{subequations}\label{eq4.11}
\begin{align}
\label{eq4.11a}
&\|e_0\|_{L^2(\mathbb{R})} \le c\bigl(\tfrac{\ell n}{2}\bigr)^{\frac{3}{2}+\epsilon-s}
\|u_0\|_{H^s_{2s}(\mathbb{R})},\\
\label{eq4.11b}
&|\mathcal{G}_n(-e_0)|^{\frac{1}{2}}
\le c \bigl(\tfrac{\ell n}{2}\bigr)^{\frac{5}{2}+\epsilon-s}
\bigl(\|u_0\|_{H^s_{2s}(\mathbb{R})}
+\|u_0\|_{H^s_{2s}(\mathbb{R})}^{\frac{5}{3}}\bigr),\\
\label{eq4.11c}
& |\mathcal{E}_{n}(e_0, 0)|^{\frac{1}{2}}\le c\bigl(\tfrac{\ell n}{2}\bigr)^{\frac{3}{2}+\epsilon-s}
\|u_0\|_{H^s_{2s}(\mathbb{R})}^{\frac{3}{2}},\\
\label{eq4.11d}
& |\mathcal{D}_n(e_0,0)|^{\frac{1}{2}}\le c\bigl(\tfrac{\ell n}{2}\bigr)^{\frac{7+2\epsilon}{4}-s}
\bigl(\|u_0\|_{H^s_{2s}(\mathbb{R})}
+\|u_0\|_{H^s_{2s}(\mathbb{R})}^{\frac{3}{2}}\bigr),\\
\nonumber
&\|\nabla_e\mathcal{D}_n(e,t)\|_{L^2([0,T]\times\mathbb{R})}
\le c\bigl(\tfrac{\ell n}{2}\bigr)^{2-s}\\
\label{eq4.11e}
&\qquad\qquad\qquad\qquad
\Bigl(\|u\|_{L^2([0,T], H^s_{2s}(\mathbb{R}))}+
\|u\|^2_{L^4([0,T], H^s_{2s}(\mathbb{R}))}\Bigr),\\
\nonumber
&\|\partial_t\nabla_e\mathcal{D}_n(e,t)\|_{L^2([0,T]\times\mathbb{R})}
\le c\bigl(\tfrac{\ell n}{2}\bigr)^{2-s}\\
\label{eq4.11f}
&\qquad\qquad\qquad\qquad
\Bigl(\|u_t\|_{L^2([0,T], H^s_{2s}(\mathbb{R}))}+
\|u_t\|^2_{L^4([0,T], H^s_{2s}(\mathbb{R}))}\Bigr).
\end{align}
\end{subequations}
In each inequality the generic constant $c>0$ is independent of $u$, $u_0$, $T>0$ and $n>0$.
\end{lemma}
\begin{proof}
(a) Since $e_0 = \mathcal{I}_n\bigl[(\mathcal{I} - \mathcal{P}_n)[u_0]\bigr]$, 
as in the proof of Corollary~\ref{lm3.4}, we obtain \eqref{eq4.11a}.

(b) We employ the Cauchy-Schwartz inequality and Lemma~\ref{lm4.2} 
(with $\varepsilon=1$) to obtain
\begin{align*}
2|\mathcal{G}_n(-e_0)| &\le c\|e_0\|^2_{L^2(\mathbb{R})} 
+ c\|\mathcal{J}_n e_0\|^2_{L^2(\mathbb{R})} + 
c\bigl(\tfrac{\ell n}{2}\bigr)^{\frac{1}{3}}\|e_0\|^{\frac{10}{3}}\\
&\le c \|(\mathcal{I}-\mathcal{I}_n)[u_0]\|^2_{H^1(\mathbb{R})}  
+c \|(\mathcal{I}-\mathcal{P}_n)[u_0]\|^2_{H^1(\mathbb{R})}  +  
c\bigl(\tfrac{\ell n}{2}\bigr)^{\frac{1}{3}}\|e_0\|^{\frac{10}{3}},
\end{align*} 
with $c>0$, depending on parameters $\alpha$, $\beta$, $\gamma$ and $\delta$ only. 
Hence, \eqref{eq4.11a} and Corollary~\ref{lm3.4} imply \eqref{eq4.11b}.

(c) From the definition of $\mathcal{E}_n(e,t)$ and Lemma~\ref{lm4.2}, we have 
\[
|\mathcal{E}_n(e_0,0)| \le 2|\delta|\|\tilde{u}_0\|_{L^\infty(\mathbb{R})}\|e_0\|^2_{L^2(\mathbb{R})}
\]
and \eqref{eq4.11c} is a consequence of Lemma~\ref{lm4.6} and \eqref{eq4.11a}.

(d) The functional $\mathcal{D}_n(e_0,0)$ is linear in $e_0$. Consequently,
\[
|\mathcal{D}_n(e_0,0)| \le \|e_0\|_{L^2(\mathbb{R})}\|\nabla_e\mathcal{D}_n(e_0,0)\|_{L^2(\mathbb{R})}
\]
and \eqref{eq4.11d} follows directly from \eqref{eq4.11a} and the proof of \eqref{eq4.11e} below.

(e) From \eqref{eq4.5c} we have 
\begin{align*}
\nabla_e\mathcal{D}_n(e,t) &= \alpha (\mathcal{I}-\mathcal{P}_n)[u]
 +\beta\mathcal{H}\bigl[(\partial_x + \mathcal{J}_n)[u]\bigr]\\
& - \gamma(\partial_{x}^2 - \mathcal{J}_n^2)[u]
 + \delta\bigl(u^2 - \mathcal{I}_n\bigl[\mathcal{P}_n[u]^2\bigr]\bigr)\\
&= \alpha E_1 + \beta E_2+\gamma E_3+\delta E_4.
\end{align*}
We bound each term separately. First of all, by Lemma~\ref{lm3.1},
\[
\|E_1\|_{L^2([0,T]\times \mathbb{R})} 
\le c\bigl(\tfrac{\ell n}{2}\bigr)^{-s}\|u\|_{L^2([0,T], H^s_{2s}(\mathbb{R}))}.
\]
Further, using the definition of $\mathcal{J}_n$, we obtain
\begin{align*}
\|E_2\|^2_{L^2(\mathbb{R})} &= 
\|\partial_x(\mathcal{I}-\mathcal{P}_n)[u]\|^2_{L^2(\mathbb{R})} + 
\tfrac{n+1}{2\ell}\bigl(|\hat{u}_{n-1}|^2 + |\hat{u}_{n}|^2\bigr)\\
&\le \|(\mathcal{I}-\mathcal{P}_n)[u]\|^2_{H^1(\mathbb{R})} + 
c\bigl(\tfrac{\ell n}{2}\bigr)\|(\mathcal{I}-\mathcal{P}_{n-2})[u]\|^2_{L^2(\mathbb{R})},
\end{align*}
so that by Lemma~\ref{lm3.1},
\[
\|E_2\|_{L^2([0,T]\times \mathbb{R})} 
\le c\bigl(\tfrac{\ell n}{2}\bigr)^{1-s}\|u\|_{L^2([0,T], H^s_{2s}(\mathbb{R}))}.
\]
Similar calculations give also
\begin{align*}
\|E_3\|_{L^2(\mathbb{R})} &\le  
\|(\mathcal{I}-\mathcal{P}_n)[u]\|_{H^2(\mathbb{R})} 
+\bigl(\tfrac{\ell n}{2}\bigr)\|(\mathcal{I}-\mathcal{P}_{n-2})[u]\|_{L^2(\mathbb{R})}
\end{align*}
and 
\[
\|E_3\|_{L^2([0,T]\times \mathbb{R})} 
\le c\bigl(\tfrac{\ell n}{2}\bigr)^{2-s}\|u\|_{L^2([0,T], H^s_{2s}(\mathbb{R}))}.
\]
Finally, 
\[
E_4 = (\mathcal{I} - \mathcal{I}_n)[u^2] + \mathcal{I}_n[u^2 - \mathcal{P}_n[u]^2] = E_{41} + E_{42}.
\]
First, we employ Lemma~\ref{lm2.2} and  Corollary~\ref{lm3.4} to obtain
\[
\|E_{41}\|_{L^2([0,T]\times \mathbb{R})}   
\le c\bigl(\tfrac{\ell n}{2}\bigr)^{\frac{3}{2}+\epsilon-s}\|u\|^2_{L^4([0,T], H^s_{2s}(\mathbb{R}))}.
\]
Second, from Lemmas~\ref{lm2.2} and \ref{lm3.2}, we infer
\begin{align*}
\|E_{42}\|_{L^2(\mathbb{R})}   
&\le c(\tfrac{\ell n}{2}\bigr) 
\|(I-\mathcal{P}_n)[u] (I+\mathcal{P}_n)[u]\|_{H^{\frac{1}{2}+\epsilon}(\mathbb{R})}\\
&\le  c(\tfrac{\ell n}{2}\bigr)\Bigl(2\|u\|_{H^{\frac{1}{2}+\epsilon}(\mathbb[R])}
+\|(I-\mathcal{P}_n)[u]\|_{H^{\frac{1}{2}+\epsilon}(\mathbb{R})}\Bigr)  
\|(I-\mathcal{P}_n)[u]\|_{H^{\frac{1}{2}+\epsilon}(\mathbb{R})}.
\end{align*}
The last bound and Lemma~\ref{lm3.1} yield
\[
\|E_{42}\|_{L^2([0,T]\times \mathbb{R})}   
\le c\bigl(\tfrac{\ell n}{2}\bigr)^{\frac{3}{2}+\epsilon-s}\|u\|^2_{L^4([0,T], H^s_{2s}(\mathbb{R}))}.
\]
Combining all our estimates together, we arrive at \eqref{eq4.11e}. 
To obtain \eqref{eq4.11f}, replace $u$ with $u_t$.
\end{proof}

Combining Theorem~\ref{lm4.5} and Lemmas~\ref{lm4.6}, \ref{lm4.7}, we obtain 
\begin{corollary}[Convergence]\label{lm4.8}
Assume $u, u_t\in L^\infty([0,T], H^s_{2s}(\mathbb{R}))$,  $s> \tfrac{11}{4}$
and $\epsilon>0$.
Then the numerical solution $\bar{u}$ satisfies
\begin{align}
\nonumber
\|u - \bar{u}\|_{L^\infty([0, T], L^2(\mathbb{R}))}
&\le c\bigl(\tfrac{\ell n}{2}\bigr)^{\frac{5}{2}+\epsilon-s}\\
\label{eq4.12}
&\Bigl(\|u_0\|_{H^s_{2s}(\mathbb{R})}
+\|u_0\|_{H^s_{2s}(\mathbb{R})}^{\frac{5}{3}}
+ \|u\|_{H^1([0,T], H^s_{2s}(\mathbb{R}))}\Bigr),
\end{align}
uniformly for large values of $n>0$, with $c>0$ that depends on the terminal time $T>0$, parameters 
$\alpha$, $\beta$, $\gamma$ and $\delta$ of the model \eqref{eq4.1} and on the regularity of the exact solution $u$
only.
\end{corollary}
\begin{proof} Note that 
\[
\|u-\bar{u}\|_{L^2(\mathbb{R})} \le 
\|e\|_{L^2(\mathbb{R})} + \|(\mathcal{I}-\mathcal{P}_n)[u]\|_{L^2(\mathbb{R})}.
\]
Hence \eqref{eq4.12} follows from Lemma~\ref{lm3.1} and the fact that under the assumption $s>\tfrac{11}{4}$, 
the numerical error $e = \tilde{u} - \bar{u}$ fells in the scope of Theorem~\ref{lm4.5}.
\end{proof}

To conclude this section, we remark that if $u, u_t\in L^\infty([0,T], H^s_{2s}(\mathbb{R}))$, for any $s>\tfrac{11}{4}$,
then, according to Corollary~\ref{lm4.8}, the convergence rate is spectral, i.e. the semi-discretization error 
$\|u - \bar{u}\|_{L^\infty([0, T], L^2(\mathbb{R}))}$ decreases faster than any inverse power of $n>0$.

\section{Implementation and simulations}\label{sec5}

The semi-discretization \eqref{eq4.1} leads to a finite dimensional system of ODEs 
whose solution is not known explicitly and itself requires an appropriate numerical treatment.
Below, we discuss briefly a suitable time-stepping algorithm and then switch to simulations.

\subsection{Implementation}
The semi-discretization \eqref{eq4.1} can be written in the form
\begin{equation}\label{eq5.1}
Y' = (\alpha J + \beta HJ^2 - \gamma J^3) Y + J F(Y) = DY+ J F(Y),\quad Y(0) = Y_0, 
\end{equation}
where, the neutral symbol $Y\in\mathbb{R}^{n+1}$, $n=2p-1$, represents either the vector 
\[
Y = (\bar{u}_0,  \ldots, \bar{u}_{2p-1}) ^T,
\] 
of the discrete MTC-Fourier coefficients or the vector of physical values 
\[
Y = (\bar{u}(x_0),\ldots,\bar{u}(x_n))^T,
\] 
computed at the nodal points $x_m$, $0\le m\le n$. The square skew-symmetric 
matrices $J, H \in\mathbb{R}^{(n+1)\times(n+1)}$ provide suitable realizations of the discrete 
operators $\mathcal{J}_n$ and $\mathcal{P}_n\mathcal{H}\mathcal{P}_n$, respectively.
The concrete form of $J$ and $H$ depends on the particular representation of $Y$.
For instance, in the MTC-Fourier (frequency) space $J$ and $H$ have simple two-by-two block structure 
with nonzero three-diagonal, respectively diagonal, blocks in the reverse block diagonal (see identities \eqref{eq3.5}), 
while both matrices are dense in the physical space. The nonlinearity $F(Y)$, 
representing $\mathcal{I}_n[\bar{u}^2]$, is given explicitly by
\[
F(Y) = \delta (\bar{u}^2(x_0),\ldots, \bar{u}^2(x_{n}))^2,
\]
in the physical space.

\paragraph{Time-stepping} The spectrum of operator $J$, computed explicitly in 
the proof of Lemma~\ref{lm4.1} (see also \cite{Wid1994}), indicates that 
\eqref{eq5.1} is stiff and hence, fully explicit time-stepping schemes cannot be unconditionally stable. 
Furthermore, since the nonlinearity $F(Y)$ is multiplied by $J$, the semi-implicit splitting-type schemes 
that separate stiff and nonstiff components of the vector field (see e.g. discussion in \cite{SPG2017}, in connection 
with the nonlinear Schr\"odinger equation) are also not plausible here. From the prospective of numerical 
stability, we are forced to use fully implicit $A$-stable algorithms. 

In our simulations, we make use of the implicit $4$-stage $8$-order Gauss-type Runge-Kutta method (IRK8 in the sequel)
\[
\begin{array}{c|c}
c& A\\
\hline
& b^T
\end{array},
\quad b,c\in\mathbb{R}^4, \quad A\in\mathbb{R}^{4\times 4},
\]
of J. Kuntzmann and  J. Butcher, (for the concrete values of the coefficients $A$, $b$ and $c$ 
see \cite[Table 7.5, p. 209]{HANOWA08}). A single IRK8 time step of length $\tau$, applied to \eqref{eq5.1}, reads
\begin{subequations}\label{eq5.2}
\begin{align}
\label{eq5.2a}
&Z_1 = \bigl(I - \tau[A\otimes D]\bigr)^{-1}\Bigl(\mathbbm{1}\otimes Y_0 + \tau [A\otimes J] G(Z_1)\Bigr),\\
\label{eq5.2b}
&Z_1 = (Y_{1,1}^T,\ldots, Y_{1,4}^T)^T,\quad 
G(Z_1) = (F(Y_{1,1})^T,\ldots, F(Y_{1,4})^T)^T,\\
\label{eq5.2c}
& Y_1 = Y_0 + \tau \sum_{i=1}^4 b_i F(Y_{1,i}),
\end{align}
\end{subequations}
where $\mathbbm{1} = (1,1,1,1)^T$ and $\otimes$ is the standard Kronecker product.
We observe that the spectrum of $A$ contains two pairs of complex conjugate eigenvalues with 
a nontrivial real part and therefore, from Lemma~\ref{lm4.1} we deduce that the spectrum of matrix 
$\bigl(I - \tau[A\otimes D]\bigr)^{-1}[A\otimes J]$ is uniformly bounded with respect to the 
space discretization parameter $n>0$. Further, the theory of Section~\ref{sec4} indicates that 
for smooth exact solutions of the Benjamin equation \eqref{eq1.1}, the semi-discrete nonlinearity $F(Y)$
is bounded uniformly in $n>0$ 
along the trajectories of \eqref{eq5.1}. Hence, the fully discrete scheme \eqref{eq5.2} is unconditionally stable.
Moreover, it follows that for a fixed $Y_0$, moderately small values of time step $\tau$ and independently of $n>0$, 
the nonlinear map, defined by the right-hand side of \eqref{eq5.2a}, is a contraction. 
As a consequence, the nonlinear equation \eqref{eq5.2a} can be solved efficiently via basic fixed point iterations. 
The observation is important from practical point of view as Newton-type iterations are prohibitively expensive 
for large values of $n>0$. 
We note also that the exact flow $\varphi_t$, generated by \eqref{eq5.1}, is symplectic. 
The IRK8 scheme is known to be symmetric and symplectic \cite{HANOWA08}, hence, the discrete 
flow of \eqref{eq5.2} preserves this property automatically.

\paragraph{Computational complexity} A single fixed point iteration, applied to \eqref{eq5.2a}, involves: solving 
linear systems with matrix $I - \tau [A\otimes D]$; the matrix-vector multiplication with matrices $D$ and $J$
and finally; computing the nonlinearity $F(Y)$. In view of the special structure of $J$ and $H$, 
in the Fourier-Christov space each matrix-vector operation requires $\mathcal{O}(n)$ flops,
while computing of $F(Y)$ involves the use of the discrete direct and inverse MTC-Fourier transforms 
(see formulas \eqref{eq3.10a} and \eqref{eq3.11}, \eqref{eq3.12}, respectively). Because of \eqref{eq3.3} and 
\eqref{eq3.10b}, both operations can be accomplished in $\mathcal{O}(n\log_2n)$ flops via the direct and 
inverse discrete Fast Fourier Transforms \cite{Boyd1990, Cri1982, IsWe2019, Wid1994, Wid1995} and the cost 
of a single iteration is $\mathcal{O}(n\log_2 n)$. As noted earlier, for any given tolerance $\varepsilon$ 
the total number of such iterations is finite and depends on the time step $\tau$ only. Hence, the overall 
complexity of a single time step of \eqref{eq5.2} is $\mathcal{O}(n\log_2 n)$.

\subsection{Simulations}

Below, we provide several simulations illustrating the accuracy of \eqref{eq4.1} in several computational 
scenarios. 

\begin{figure}[t]
\begin{center}
\begin{minipage}{0.45\textwidth}
\includegraphics[scale=0.9]{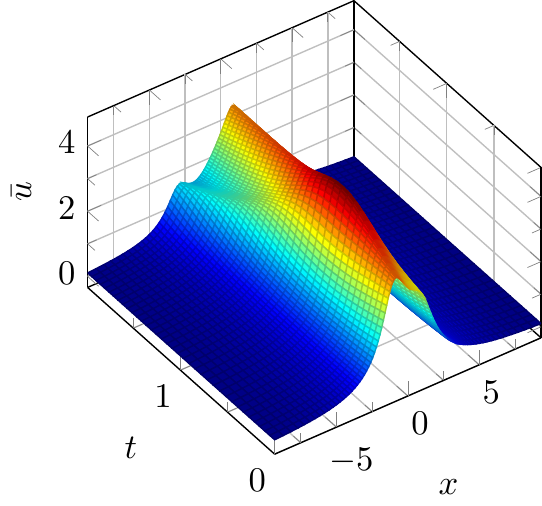}
\includegraphics[scale=0.9]{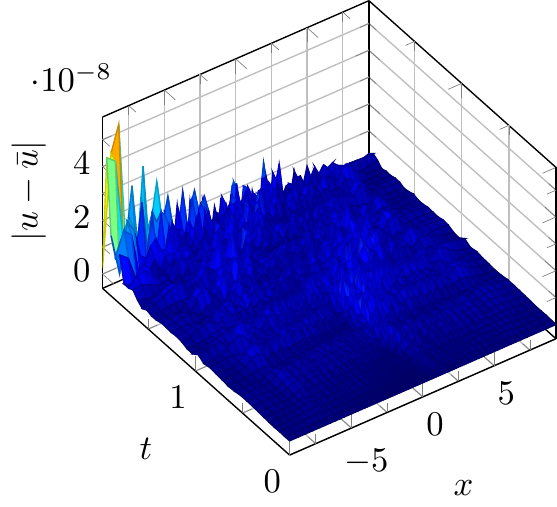}
\end{minipage}\hspace{0.5cm}
\begin{minipage}{0.45\textwidth}
\includegraphics[scale=0.9]{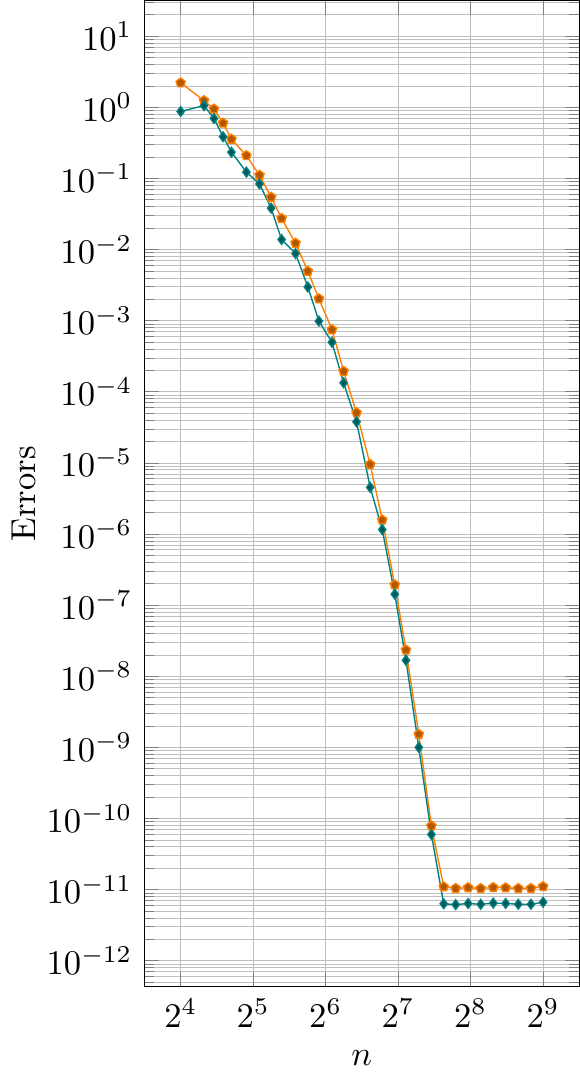}
\end{minipage}
\end{center}
\caption{The left diagrams (top to bottom): the numerical solution of the Benjamin equation 
\eqref{eq1.1} $\bar{u}$ and the pointwise error $|u-\bar{u}|$, with $n=2^7-1$.  
The right diagram: $\|u-\bar{u}\|_{L^\infty([0,T];L^2(\mathbb{R}))}$ (orange, pentagon), 
$\|u-\bar{u}\|_{L^\infty([0,T]\times \mathbb{R})}$ (teal, diamond). }\label{fig5.1}
\end{figure}

\subsubsection{Slowly decreasing solutions}
We begin with the generic situation whe\-re, due to the nature of the Fourier symbol 
in the linear part of \eqref{eq1.1}, solutions decay at most algebraically.

\paragraph{Example 1} First, we simulate \eqref{eq1.1} in time interval $[0,2]$, with 
$\alpha = \beta = \gamma = \delta = 1$. Since for these values of the model parameters, analytic formulas 
for solutions are not available, we augment \eqref{eq1.1a} with a source term $f(x,t)$. The latter is chosen 
so that the exact solution reads
\begin{align*}
&u(x,t) = \sum_{k=1}^3 \tfrac{r_k}{a_k^2+(x-x_{k,0}-c_kt)^2},\\
&r_1 = 2,\quad r_2=1,\quad r_3 = 3,\quad
a_1 = 1,\quad a_2 = 1,\quad a_3 = 2,\\
&c_1 = 1,\quad c_2 = -2,\quad c_3 = 0,\quad
x_{1,0} = -1,\quad x_{2,0} = 1,\quad x_{3,0} = 0.
\end{align*} 
Note that $u(x,t)$ is smooth (in fact $u\in  H^s(\mathbb{R})$, $s\in\mathbb{R}$), 
but has a polynomial decay rate at infinity ($u = \mathcal{O}(|x|^{-2})$ at $|x|\to\infty$). 
In view of this fact, accurate approximation of such functions with the aid of standard 
trigonometric basis requires huge number of spatial grid points.
Nevertheless, straightforward calculations show that the quantity $\mathcal{F}\mathcal{P}_+[u]$ 
is smooth and decreases exponentially in the positive half line $\mathbb{R}_+$. 
Hence, $u$ falls in the scope of the theory presented in Sections~\ref{sec3} and \ref{sec4} 
and we expect rapid error decay already for moderate values of $n>0$.

The numerical results, for $2^4-1\le n\le 2^9-1$, $\ell=2^3$ and $\tau = 2\cdot10^{-2}$, 
are plotted in the right diagram of Fig.~\ref{fig5.1}. Both $\|\cdot\|_{L^\infty([0,T], L^2(\mathbb{R}))}$ 
and $\|\cdot\|_{L^\infty([0,T]\times \mathbb{R})}$ errors decrease spectrally (note that both curves are concave) 
as $n$ increases. For $n>2^7$, the numerical errors settle near $10^{-11}$. This is a consequence of the 
inexact time-stepping procedure employed in our calculations. Simulations, not reported here, indicate that 
for $n>2^7$ the error can be further reduced by choosing smaller time integration steps.

To illustrate the quality of the approximation, we plot the numerical solution $\bar{u}$ and 
the associated pointwise error $|u-\bar{u}|$, obtained with $n=2^7-1$, in the  two left diagrams 
of Fig.~\ref{fig5.1}.  It is clearly visible that the pointwise error does not exceed the magnitude 
of $5\cdot 10^{-8}$ uniformly in the computational domain.

\begin{figure}[t]
\begin{center}
\begin{minipage}{0.45\textwidth}
\includegraphics[scale=0.9]{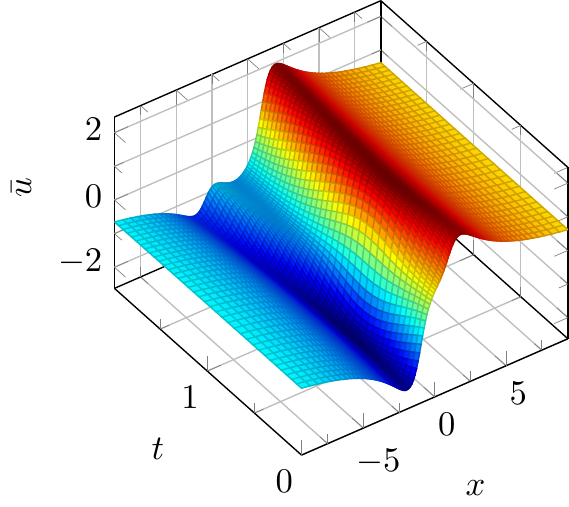}
\includegraphics[scale=0.9]{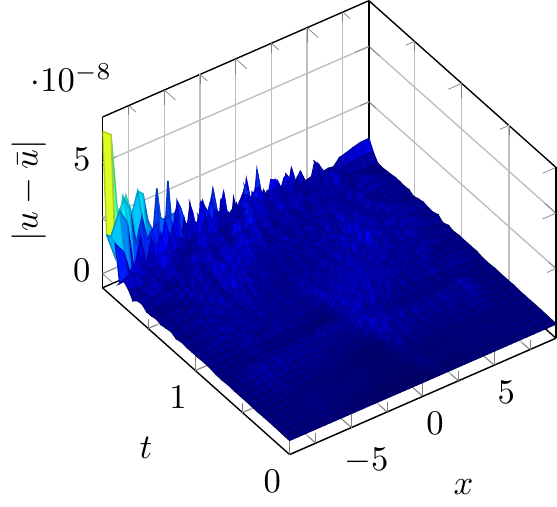}
\end{minipage}\hspace{0.5cm}
\begin{minipage}{0.45\textwidth}
\includegraphics[scale=0.9]{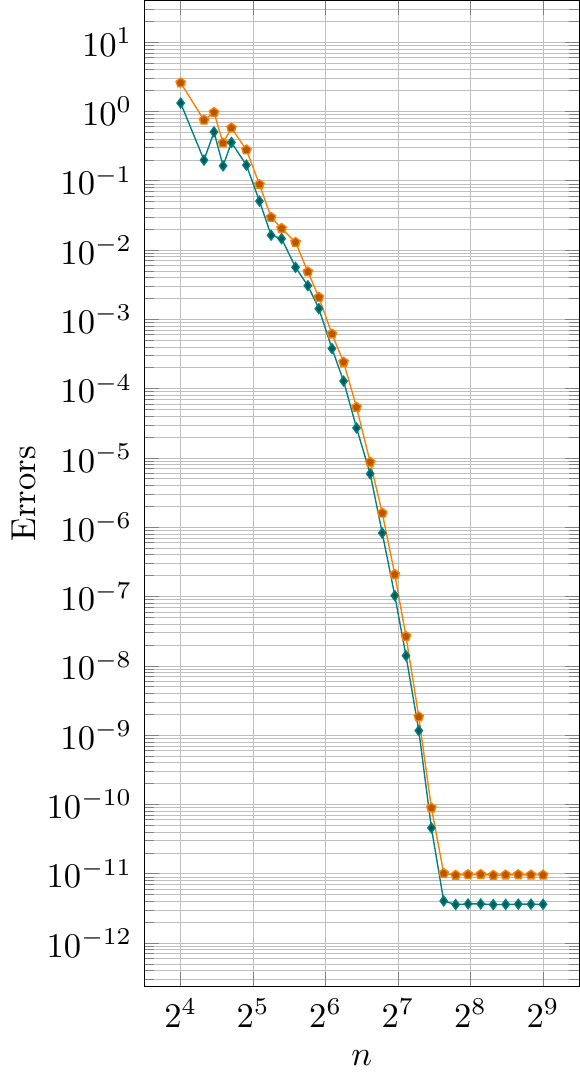}
\end{minipage}
\end{center}
\caption{The left diagrams (top to bottom): the numerical solution of the Benjamin equation 
\eqref{eq1.1} $\bar{u}$ and the pointwise error $|u-\bar{u}|$, with $n=2^7-1$.  
The right diagram: $\|u-\bar{u}\|_{L^\infty([0,T];L^2(\mathbb{R}))}$ (orange, pentagon), 
$\|u-\bar{u}\|_{L^\infty([0,T]\times \mathbb{R})}$ (teal, diamond). }\label{fig5.2}
\end{figure}

\paragraph{Example 2} In our second simulation,  we keep the numerical parameters of Example 1
unchanged, but make use of another source term which gives the following exact solution
\[
u(x,t) = \sum_{k=1}^3 \tfrac{r_k(x-x_{k,0}-c_kt)}{a_k^2+(x-x_{k,0}-c_kt)^2}.
\] 
In this settings $u(x,t) = \mathcal{O}(|x|^{-1})$, as $|x|\to\infty$. Nevertheless, the truncated 
Fourier image $\mathcal{F}\mathcal{P}_+[u]$ has exactly the same qualitative features 
as in Example 1 and the resulting convergence rate is spectral (see the left diagram in 
Fig~\ref{fig5.2}). In the particular case of $n=2^7-1$, the numerical solution and the pointwise 
error are shown in the top- and bottom-left diagrams of Fig.~\ref{fig5.2}, respectively. The observed 
behavior is very much alike to the one, reported in Example 1.

\begin{figure}[t]
\begin{center}
\begin{minipage}{0.45\textwidth}
\includegraphics[scale=0.9]{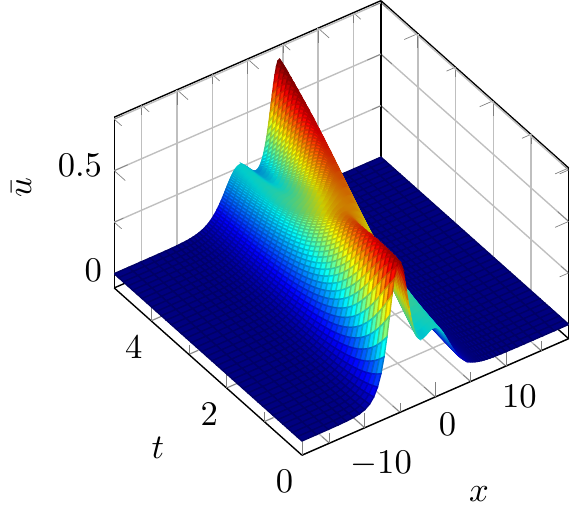}
\includegraphics[scale=0.9]{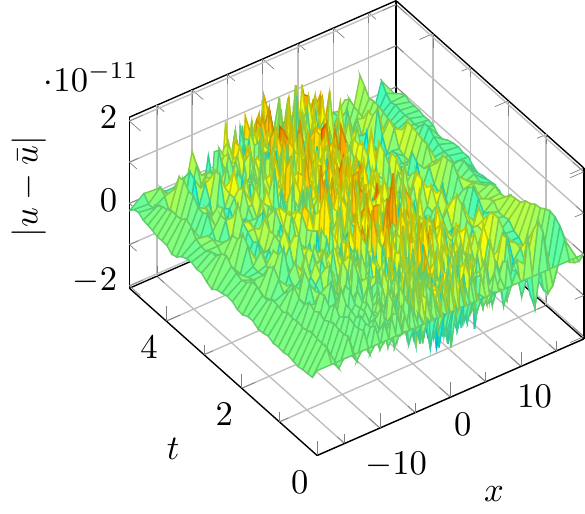}
\end{minipage}\hspace{0.5cm}
\begin{minipage}{0.45\textwidth}
\includegraphics[scale=0.9]{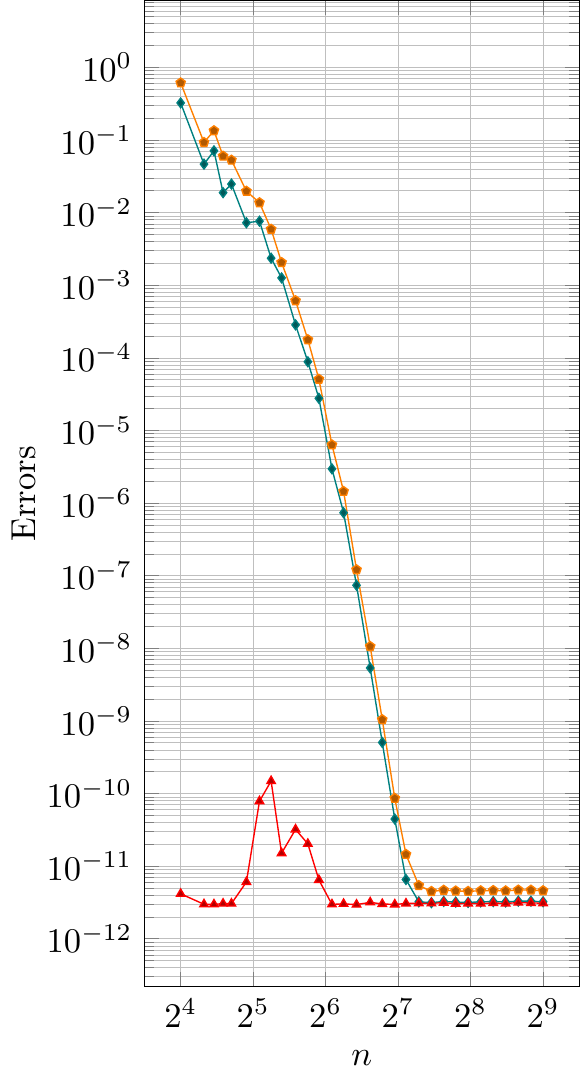}
\end{minipage}
\end{center}
\caption{ The left diagrams (top to bottom): the numerical solution of the Benjamin equation 
\eqref{eq1.1} $\bar{u}$ and the pointwise error $|u-\bar{u}|$, with $n=2^7-1$.  
The right diagram: $\|u-\bar{u}\|_{L^\infty([0,T];L^2(\mathbb{R}))}$ (orange, pentagon), 
$\|u-\bar{u}\|_{L^\infty([0,T]\times \mathbb{R})}$ (teal, diamond). }\label{fig5.3}
\end{figure}

\subsubsection{The Korteweg-de Vries scenario}
The Benjamin equation \eqref{eq1.1} contains two special case 
$\gamma = 0$ and $\beta = 0$, which are of independent interest. 
The first one corresponds to the Benjamin-Ono equation, and is not considered here. 
In the second case, we have the classical Korteweg-de Vries (KdV) equation. 
The latter is known to be completely integrable and possesses a large number of special solutions.
For instance, when
\[
\alpha=\beta=0,\quad \gamma = -1,\quad \delta=-3,
\]
the inverse scattering transform yields the so called $N$-solitons (see e.g. \cite{Abl2004})
\begin{subequations}\label{eq5.3}
\begin{align}
\label{eq5.3a}
&u(x,t) = - 2\partial_{xx} \ln \det(I + A(x,t)),\\
\label{eq5.3b}
&A(x,t) = \bigl(b_ie^{8 \lambda^3_i t} 
\tfrac{e^{-\lambda_ix-\lambda_j x}}{\lambda_i+\lambda_j}, 1\le i,j\le N\bigr),\\
\label{eq5.3c}
&\lambda_i = \tfrac{1}{2}\sqrt{v_i}\quad
b_i =   2\lambda_ie^{2\phi_i\lambda_i},\quad 1\le i\le N,
\end{align}
\end{subequations}
which describe evolution of $N$ traveling waves, whose velocities and the phases
are controlled by $v_i$ and $\phi_i$, respectively. Directly from \eqref{eq5.3}, it follows that 
$N$-solitons are smooth and decay exponentially to zero as $|x|$ increases.
Hence, such solutions fall in the scope of the theory developed in Sections~\ref{sec3}~and~\ref{sec4}.

\paragraph{Example 3} To illustrate the above statement, in \eqref{eq5.3} we let 
\[
v_1 = \tfrac{3}{2},\quad v_2= \tfrac{1}{2},\quad \phi_1 = -3,\quad \phi_2 = 0,
\]
choose $u_0$ according to \eqref{eq5.3}, take $\ell=2^3$, $\tau = 10^{-2}$ and integrate 
\eqref{eq4.1} numerically in time interval $[0,5]$. The results of simulations (see in Fig.~\ref{fig5.3}) 
are qualitatively similar to those obtained in Examples 1 and 2.
In particular, the plots of $\|\cdot\|_{L^\infty([0,T], L^2(\mathbb{R}))}$ and 
$\|\cdot\|_{L^\infty([0,T]\times \mathbb{R})}$  errors indicate that the convergence rate 
is spectral. Note however that in the bottom-left diagram of Fig.~\ref{fig5.3} the pointwise
error is smaller than in the two previous Examples. This is connected with the exponential 
decay of the $2$-soliton at infinity (its accurate spatial resolution requires fewer grid points 
than in Examples 1 and 2).

By construction, the scheme \eqref{eq4.1} is conservative and the semi-discrete Hamiltonian 
$\mathcal{G}_n(\bar{u})$ remains constant along the exact trajectories of \eqref{eq4.1}.
In order to test the conservation properties of the fully discrete scheme, in the right diagram 
of Fig.~\ref{fig5.3}, we added the plot of the quantity $\max_{t\in[0,T]}|\mathcal{G}_n(\bar{u}_0) 
- \mathcal{G}_n(\bar{u})|$, measuring the largest deviation in the Hamiltonian. We observe 
that the deviation remains several orders of magnitude smaller than 
either of the $\|\cdot\|_{L^\infty([0,T], L^2(\mathbb{R}))}$ and 
$\|\cdot\|_{L^\infty([0,T]\times \mathbb{R})}$  errors, until the latter settle near $10^{-11}$.

\begin{figure}[t]
\begin{center}
\begin{minipage}{0.45\textwidth}
\includegraphics[scale=0.9]{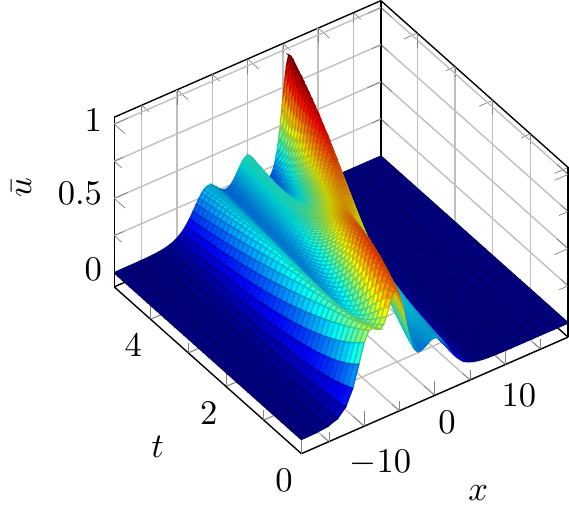}
\includegraphics[scale=0.9]{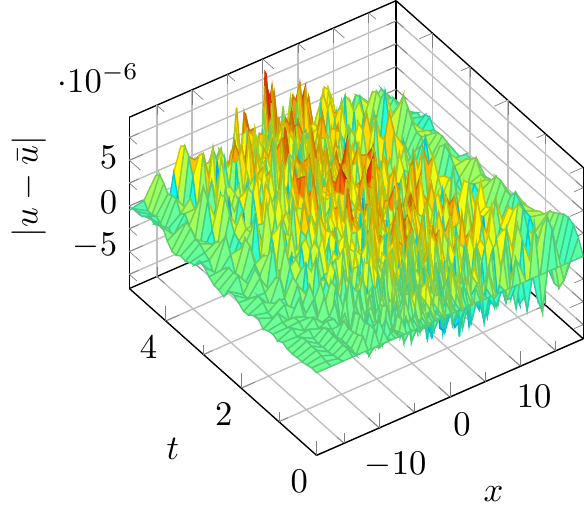}
\end{minipage}\hspace{0.5cm}
\begin{minipage}{0.45\textwidth}
\includegraphics[scale=0.9]{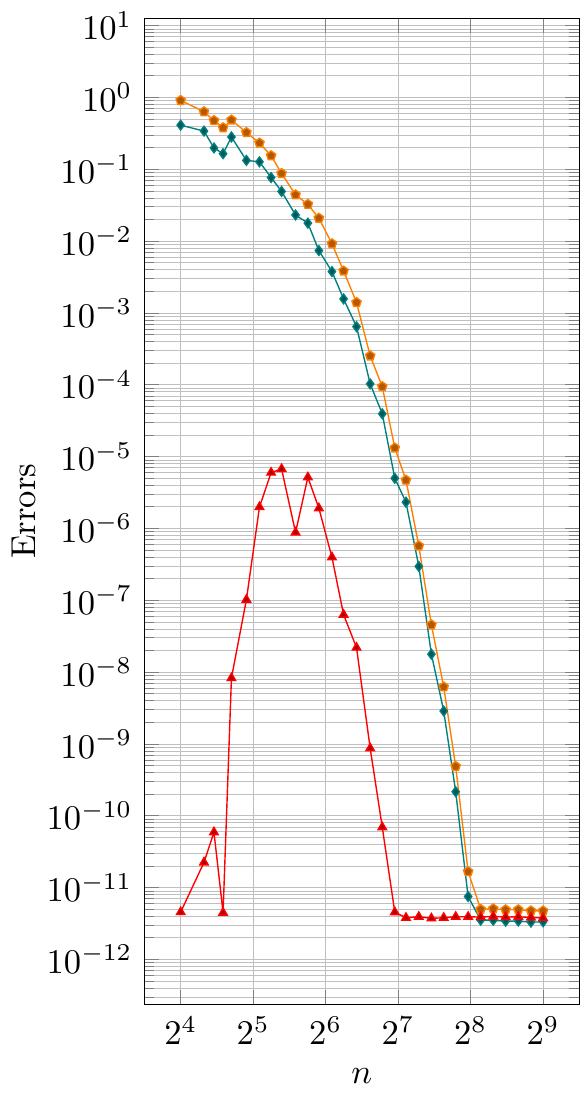}
\end{minipage}
\end{center}
\caption{The left diagrams (top to bottom): the numerical solution of the Benjamin equation 
\eqref{eq1.1} $\bar{u}$ and the pointwise error $|u-\bar{u}|$, with $n=2^7-1$.  
The right diagram: $\|u-\bar{u}\|_{L^\infty([0,T];L^2(\mathbb{R}))}$ (orange, pentagon), 
$\|u-\bar{u}\|_{L^\infty([0,T]\times \mathbb{R})}$ (teal, diamond). }\label{fig5.4}
\end{figure}

\paragraph{Example 4} We repeat calculations of Example 3, but this time with
\[
v_1 = 1,\quad v_2= 1,\quad v_3 = \tfrac{1}{2},\quad \phi_1 = -4,\quad \phi_2 = -2,\quad \phi_3 = 0.
\]
This scenario describes an elastic collision of three traveling waves, see the top-left diagram in Fig.~\ref{fig5.4}. 
The exact $3$-soliton has exactly the same qualitative features as the $2$-soliton of Example 3, with the 
exception that now the exponential decay rate is slightly slower. This manifests in larger numerical 
errors, see the bottom-left diagram in Fig.~\ref{fig5.4}.

\begin{figure}[t]
\begin{center}
\includegraphics{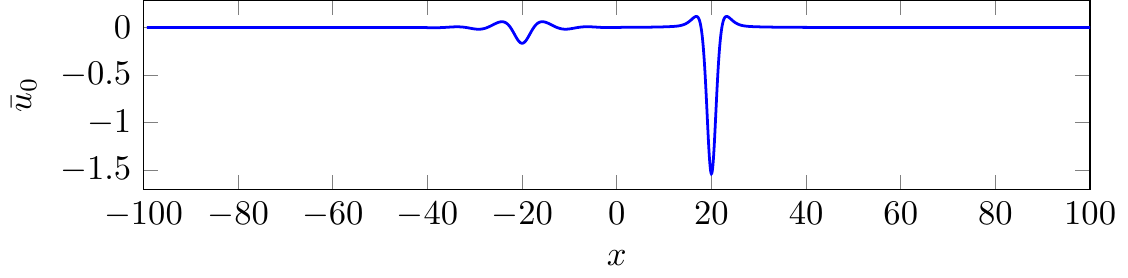}
\includegraphics{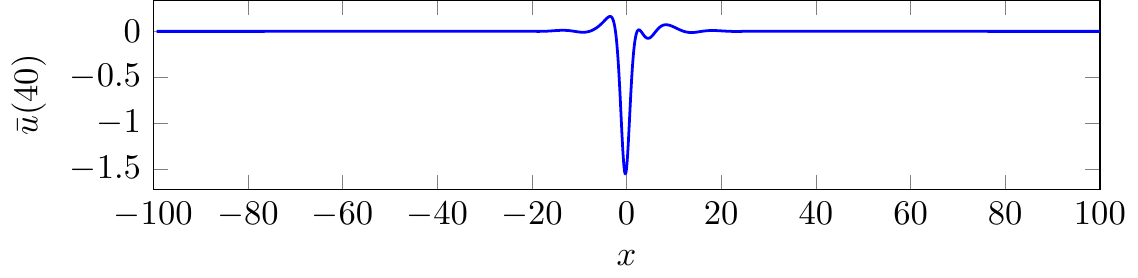}
\includegraphics{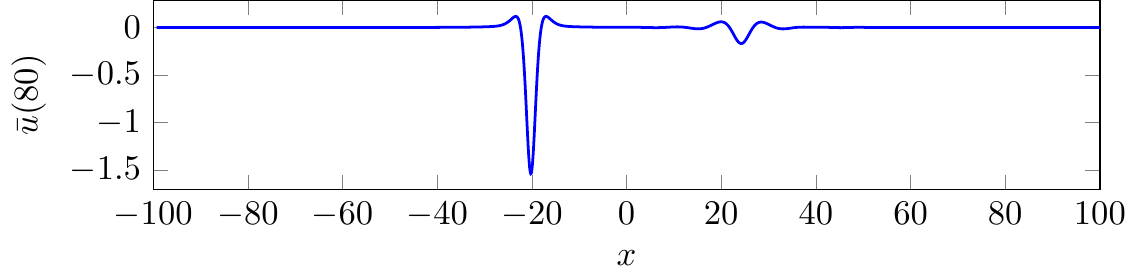}
\includegraphics{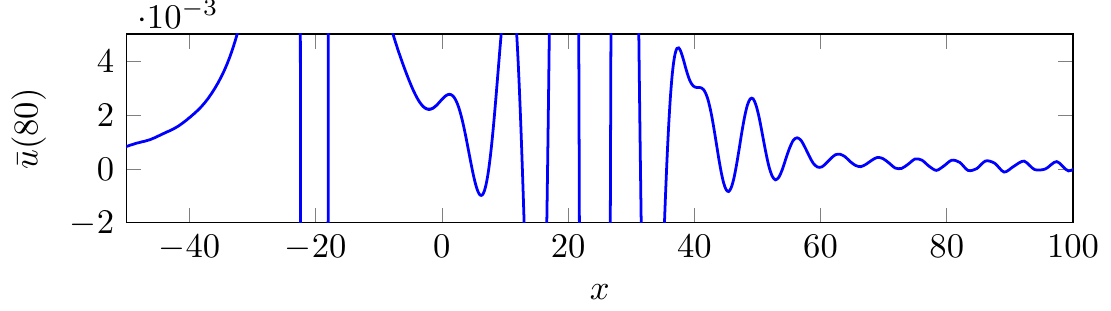}
\end{center}
\caption{Inelastic collision of two traveling waves, Example 5.}\label{fig5.5}
\end{figure}

\subsubsection{Traveling waves}
In our last two simulations, we model an interaction of traveling waves. In the context of the
Benjamin equation \eqref{eq1.1}, the traveling wave solutions are given by 
$u(x,t) = v_\sigma(x-ct)$, where $v_\sigma$ satisfies
\begin{subequations}\label{eq5.4}
\begin{align}
\label{eq5.4a}
&v_\sigma - 2\sigma \sqrt{\tfrac{\gamma}{\alpha-c}}\mathcal{H}[\partial_x v_\sigma] 
- \tfrac{\gamma}{\alpha-c}\partial_{xx}v_\sigma 
+ \tfrac{\delta}{\alpha-c}v^2_\sigma= 0,\quad x\in\mathbb{R},\\
\label{eq5.4b}
&\sigma = \tfrac{\beta}{2\sqrt{\gamma(\alpha-c)}},\quad \gamma,\delta,\nu>0,\quad c<\alpha.
\end{align}
\end{subequations}
For a rigorous treatment of \eqref{eq5.4}, see 
\cite{BENJ1996, AlBoRe1999,Pava1999, KaBo2000,CaAk2003,DoDuMi2015, DoDuMi2016} and references therein.

The exact solutions to \eqref{eq5.4}, apart from the trivial case of $\alpha = 0$, are not available.
In our simulations, the \emph{even} traveling waves are constructed numerically. We employ a simple continuation 
scheme, which works as follows: for a given $\alpha$, $\beta$, $\gamma$, $\delta$ and $c$, 
that satisfy \eqref{eq5.4b} and $0\le \sigma<1$; (i) we let 
\[
\bar{v}_0 = \mathcal{I}_n[v_0],\quad v_0(x) = - \tfrac{3(\alpha-c)}{2\delta}
\sech\bigl(\sqrt{\tfrac{\alpha-c}{4\gamma}}x\bigr)^2;
\]  
(ii) introduce a continuation grid $0<\sigma_1<\ldots<\sigma_N= \sigma$ and 
(iii) apply simplified Newton's iterations to the sequence of the discrete nonlinear problems
\begin{equation}\label{eq5.5}
\bar{v}_{\sigma_j} + 2\sigma_{j} \sqrt{\tfrac{\gamma}{\alpha-c}}\mathcal{H}[\mathcal{J}_n \bar{v}_{\sigma_j}] 
- \tfrac{\gamma}{\alpha-c}\mathcal{J}^2_n\bar{v}_{\sigma _j}
+ \tfrac{\delta}{\alpha-c}\mathcal{I}_n[\bar{v}^2_{\sigma_j}]= 0,\quad 1\le j\le N,
\end{equation}
where for each $j$, $\bar{v}_{\sigma_j}$ is restricted to be even. The iterations terminate
when the $L^2(\mathbb{R})$-norm of the defect in \eqref{eq5.5} drops below the accuracy 
threshold of $\varepsilon_n = 10^{-12} \sqrt{\tfrac{2(1-\sigma)}{n}}$.
A comprehensive analysis of \eqref{eq5.5} falls outside the scope of the present paper, we mention 
only that in all our simulations the simplified Newton's process converges rapidly to the discrete 
solutions $\bar{v}_{\sigma_j}$ but, as observed by many authors, the number of iterations 
increases when $\sigma$ approaches its upper bound~of~$1$.

\paragraph{Example 5} We let $n=2^{12}-1$, $\ell=2^3$, $\alpha=\gamma=\delta=1$,
$c_1 = \tfrac{1}{2}$, $c_2 = -\tfrac{1}{2}$, $\sigma_1 = 0.95$,
$\beta = \sigma_1\sqrt{4\gamma(\alpha-c_1)}$ and $\sigma_2 = \tfrac{\beta}{\sqrt{4\gamma(\alpha-c_2)}}$.
As an initial condition, we take the sum of two translated traveling waves 
\[
\bar{u}_0(x) = \bar{v}_{\sigma_1}(x+20) + \bar{v}_{\sigma_2}(x-20)
\]
and integrate \eqref{eq4.1} numerically in time interval $[0,80]$. With this settings, the solution  
describes a collision of two traveling waves moving towards each other. The collision occurs near 
$t = 40$, past that time the waves continue to move in the opposite directions. 
The initial profile of the numerical solution and its profiles near the collision time and at the terminal time 
are shown in the top three diagrams of Fig.~\ref{fig5.5}. As observed in \cite{KaBo2000, DoDuMi2015}, 
an interaction of the Benjamin traveling waves is inelastic --- after collisions, numerical solutions 
develop a persistent small amplitude oscillating tail.  In agreement with these observations, 
the latter is clearly visible in the bottom diagram of Fig.~\ref{fig5.5}, where the magnified view 
of the terminal profile is presented.

\begin{figure}[t]
\begin{center}
\includegraphics{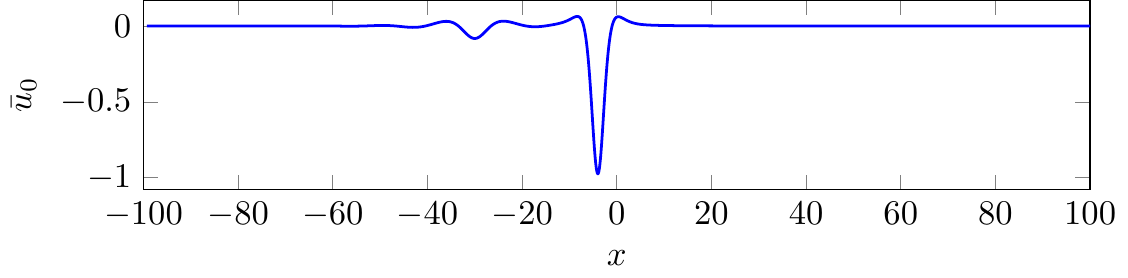}
\includegraphics{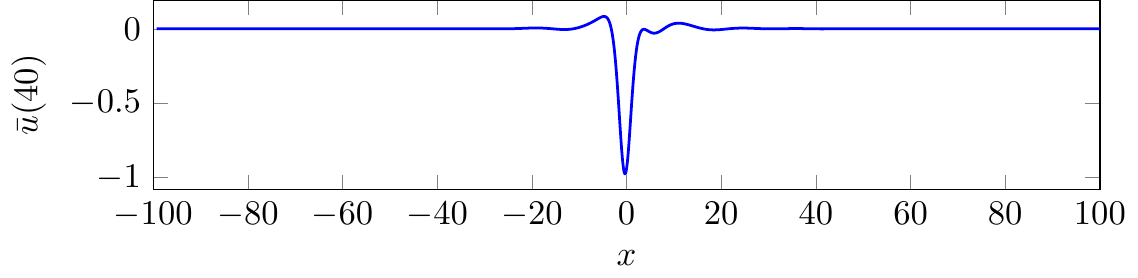}
\includegraphics{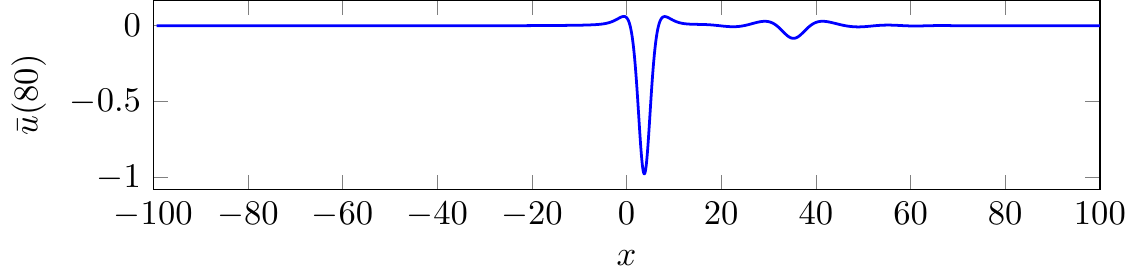}
\includegraphics{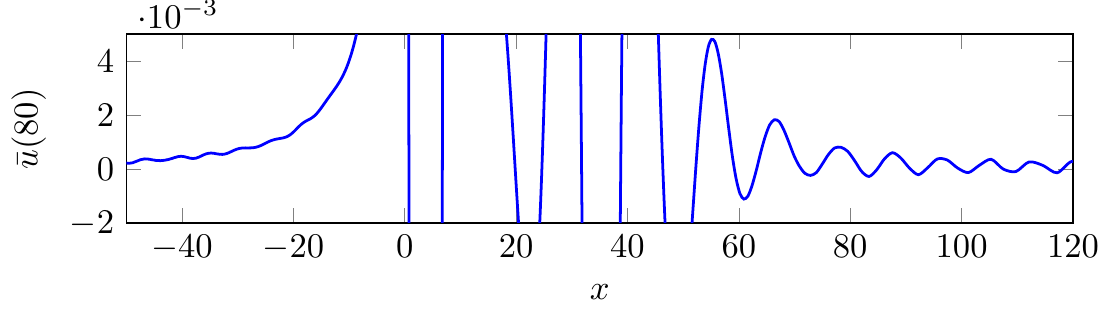}
\end{center}
\caption{Inelastic collision of two traveling waves, Example 6.}\label{fig5.6}
\end{figure}

\paragraph{Example 6}  In our last example, we use $n=2^{12}-1$, $\ell=2^3$, $\alpha=\gamma=\delta=1$,
$c_1 = \tfrac{3}{4}$, $c_2 = \tfrac{1}{10}$, $\sigma_1 = 0.95$,
$\beta= \sigma_1\sqrt{4\gamma(\alpha-c_1)}$, $\sigma_2 = \tfrac{\beta}{\sqrt{4\gamma(\alpha-c_2)}}$,
\[
\bar{u}_0(x) = \bar{v}_{\sigma_1}(x+30) + \bar{v}_{\sigma_2}(x+4)
\]
and $[0, 80]$ for the time integration interval. The scenario describes propagation of two traveling waves
moving in the same direction and colliding near $t = 40$. 
The numerical results are shown in Fig.~\ref{fig5.6}, where as before, the top three diagrams 
contain the solution profiles at the initial, near collision and the terminal times, while 
the bottom diagram contains a magnified view of the solution at terminal time $t=80$. 
Once again, the small 
dispersive tails  (of the amplitudes $\approx10^{-4}$ before the slow wave 
and $\approx10^{-3}$ after the fast wave) are clearly visible.
\section{Conclusion}\label{sec6}

In the paper, we study in detail approximation properties of the Malmquist-Takenaka-Christov (MTC) 
system. Theoretical analysis of Sections~\ref{sec3} indicates that the MTC approximations 
converge rapidly, provided Fourier images of functions, being approximated, are regular away 
from the origin and decay rapidly at infinity. The latter situation is generic for solutions of semi- and 
quasi-linear equations containing Fourier multipliers, whose symbols are irregular at the origin. 
Typical examples are models involving Hilbert/Riesz transforms (e.g. the Benjamin and the 
Benjamin-Ono equations), fractional derivatives (e.g. fractional dispersion or diffusion), e.t.c. 
We believe that for such problems the spectral/collocation MTC schemes have clear theoretical 
advantage over conventional trigonometric-Fourier, Hermite or algebraically mapped Chebyshev 
spectral approximations. 

It is worth mentioning, that unlike earlier approximation results \cite{Boyd1990, Wid1995}, we 
derive MTC error bounds directly in the functional settings of $H^s_r(\mathbb{R})$. As far as we are aware,
these results are new and can be used directly in a theoretical analysis of spectral/collocation MTC schemes.
As a concrete application, in Sections~\ref{sec4} and \ref{sec5} we provide a comprehensive treatment of 
the nonlinear Benjamin equation. In particular, we demonstrate that the MTC collocation scheme converges 
rapidly and admits an efficient implementation, comparable to the best spectral Fourier and hybrid spectral 
Fourier/finite-element methods published in the literature up to the date.

Though in the paper we mainly deal with the analysis of the MTC system and its applications, 
Appendix~\ref{sec7} contains some extensions of recent results in the theory of weighted function spaces, 
which are of independent interest.

\appendix\section{Proof of \eqref{eq2.5}}\label{sec7}

In this section, \eqref{eq2.5} is derived as as a consequence of a more general result, concerning complex 
interpolation of weighted Bessel potential spaces in $\mathbb{R}_+$. In our proof, we combine 
the notion of one-sided $A_{\pm,p}$ classes of \cite{Saw1986} with the localization ideas of \cite{Ry2001}.

\subsection{$A_{p,\pm}^{\loc}(\Omega)$ weights}
Let $\Omega$ be an open subset of $\mathbb{R}$ and $w\in L^{1,\loc}(\Omega)$ be a.e. positive 
function (weight) in $\Omega$. With $w$ and parameter $1<p<\infty$, we associate new weight 
$\bar{w}_p = w^{-\frac{p'}{p}}$, fix some $t>0$ and define
\begin{subequations}\label{eq7.1}
\begin{align}
\label{eq7.1a}
&[w]_{p,t}^+  = \sup_{[a,b]\subset\Omega,0<|b-a|<t,\,x\in(a,b)} 
\frac{w^{\frac{1}{p}}([a,x]) \bar{w}_p^{\frac{1}{p'}}([x,b])}{|b-a|},\\
\label{eq7.1b}
&[w]_{p,t}^-  = \sup_{[a,b]\subset\Omega,0<|b-a|<t,\,x\in(a,b)} 
\frac{\bar{w}_p^{\frac{1}{p'}}([a,x]) w^{\frac{1}{p}}([x,b])}{|b-a|}.
\end{align}
\end{subequations}
The $A_{p,\pm}^{\loc}(\Omega)$, $1<p<\infty$, class consists of all a.e. positive locally integrable 
functions $w$ in $\Omega$, such that for some $t<0$, the quantity $[w]^{\pm}_{p,t}$ is finite.\footnote{
$A_{p,\pm}^{\loc}(\Omega)$ class is the local version of $A_{\pm,p}$ weights of E. Sawyer, 
introduced in connection with one-sided Hardy-Littlewood maximal functions in \cite{MaRe1993, Saw1986}. 
The idea of localizing the $A_p$ condition of B. Muckenhoupt is due to V.S. Rychkov, see \cite{Ry2001}.} 
The following result can be viewed as an analogue  of \cite[Lemma~1.2]{Ry2001}.
\begin{lemma}\label{lm7.1}
The classes $A_{p,\pm}^{\loc}(\Omega)$ are well defined, i.e.
independent of a particular choice of the cut-off parameter $t>0$. 
\end{lemma}
\begin{proof}
(a) Assume for a fixed $t_0>0$, $[w]_{p,t_0}^+ = [\bar{w}_p]_{p',t_0}^-$ is finite. 
Then for any $[a,b]\subset \Omega$, with 
$0<|b-a|<t_0$ and any $x\in (a,b)$, formula \eqref{eq7.1a}, combined with H\"older's inequality, implies 
\[
\mu^p((x,b)) w([a,x]) \le ([w]_{p,t_0}^+)^p \mu^p([a,b]) w([x,b]). 
\]
This allow us to conclude that 
\begin{subequations}\label{eq7.2}
\begin{equation}\label{eq7.2a}
\Bigl(\tfrac{|x-b|}{|b-a|}\Bigr)^p \le \bigl(([w]_{p,t_0}^+)^p+1\bigr)\tfrac{w([x,b])}{w([a,b])}.
\end{equation}
In complete analogy, formula \eqref{eq7.1b} implies
\begin{equation}\label{eq7.2b}
\Bigl(\tfrac{|x-a|}{|b-a|}\Bigr)^{p'} 
\le \bigl(([\bar{w}_p]_{p',t_0}^-)^{p'}+1\bigr)\tfrac{\bar{w}_p([a,x])}{\bar{w}_p([a,b])}.
\end{equation}
\end{subequations}

(b) Consider now an interval $[a',b']\subset \Omega$, $0<|b'-a'|<2t_0$ and chose
$a'<a<x_1<b<b'$ so that $|a-a'| = |x_1-a| = |b-x_1| = |b'-b|$. 
By virtue of \eqref{eq7.2} and \eqref{eq7.1a}, we have the estimate
\[
w^{\frac{1}{p}}([a',x_1]) \bar{w}_p^{\frac{1}{p'}}([x_1,b']) 
\le c' \mu([a',b']),
\]
with $c' = 4 [w]_{p,t_0}^+\bigl(([w]_{p,t_0}^+)^p+1\bigr)^{\frac{1}{p}} 
\bigl(([w]_{p,t_0}^+)^{p'}+1\bigr)^{\frac{1}{p'}}$.
When $x\in(a',b')$ is arbitrary, the last inequality yields either
\[
w^{\frac{1}{p}}([a',x]) \bar{w}_p^{\frac{1}{p'}}([x,b']) \le
w^{\frac{1}{p}}([a',x]) \bar{w}_p^{\frac{1}{p'}}([x,x_1])  + 
w^{\frac{1}{p}}([a',x_1]) \bar{w}_p^{\frac{1}{p'}}([x_1,b']),  
\]
for $a'<x<x_1$, or 
\[
w^{\frac{1}{p}}([a',x]) \bar{w}_p^{\frac{1}{p'}}([x,b']) \le 
w^{\frac{1}{p}}([x_1,x]) \bar{w}_p^{\frac{1}{p'}}([x_1,b'])  + 
w^{\frac{1}{p}}([a',x_1]) \bar{w}_p^{\frac{1}{p'}}([x_1,b']), 
\]
for $x_1\le x < b'$. In both cases, we obtain
\[
w^{\frac{1}{p}}([a',x]) \bar{w}_p^{\frac{1}{p'}}([x,b']) \le
c'' \mu([a',b']),
\]
where $c''=  [w]_{p,t_0}^+\biggl(1 + 4\bigl(([w]_{p,t_0}^+)^p+1\bigr)^{\frac{1}{p}} 
\bigl(([w]_{p,t_0}^+)^{p'}+1\bigr)^{\frac{1}{p'}}\biggr)$. 
Hence, $[w]_{p,2t_0}^+\le c'' <\infty$.
\end{proof}

Properties of $A_{p,\pm}^\loc(\Omega)$ and $A_{p,\pm}$ weights are quite similar. 
In particular, the one-sided local Hardy-Littlewood maximal functions 
\begin{subequations}\label{eq7.3}
\begin{align}
\label{eq7.3a}
\mathcal{M}^+_{t}[f](x) = \sup_{[x,b]\subset\Omega, 0<|b-x|<t} \tfrac{1}{|b-x|} \int_x^b |f| dx,\\
\label{eq7.3b}
\mathcal{M}^-_{t}[f](x) = \sup_{[a,x]\subset\Omega, 0<|x-a|<t} \tfrac{1}{|x-a|} \int_a^x |f| dx,
\end{align}
are bounded from $L^p_w(\Omega)$, $1<p<\infty$, into itself, i.e
\begin{equation}\label{eq7.3c}
\|\mathcal{M}^\pm_{t}[f]\|_{L^p_w(\Omega)} \le c_{t,w} \|f\|_{L^p_w(\Omega)}, \quad 
1<p<\infty,\quad t>0,
\end{equation}
\end{subequations}
 if and only if $w\in A_{p,\pm}^\loc(\Omega)$.
The claim follows e.g. from the verbatim repetition of the arguments presented in \cite{MaRe1993, Saw1986}.\footnote{ 
In the context of local Muckenhoupt classes $A_p^\loc$, such results are obtained via local extensions 
of a weight $w\in A_p^\loc$ to $\bar{w}\in A_p$, see Lemma~1.1 in \cite{Ry2001}. 
This approach is not plausible in the one-sided settings, due to the asymmetric nature of \eqref{eq7.1}, 
the adjoint weight $\bar{w}_{p}$ might have non-integrable singularities at the boundary points of $\Omega$.}

\subsection{Bessel potential spaces with $A_{p,\pm}^{\loc}(\mathbb{R}_\pm)$ weights}
The Bessel potential spaces in $\mathbb{R}_\pm$ (see \cite{BaPaPoSh2014}) are defined as the images of the weighted 
Lebesgue spaces $L^{p}_w(\mathbb{R}_\pm)$ under the action of one-sided Bessel fractional integrals
$\mathcal{J}_\mp^{s}[\cdot]$, i.e. $L^{p,s}_w(\mathbb{R}_\pm) = 
\mathcal{J}_{\mp}^{s}[L^p_w(\mathbb{R}_\pm)]$, $s\ge 0$, (see Section~\ref{sec2} also).

For a fixed $t>0$ and $\varphi^\pm\in C^{\infty}_0(\mathbb{R})$, with
$\supp\varphi^\pm \subset \mathbb{R}_{\pm}\cap (-t,t)$, we define 
\[
\mathcal{T}_{\varphi^\mp,t}[f](x) = \varphi^\mp \ast f,\quad x\in\mathbb{R}_+.
\]
\begin{lemma}\label{lm7.2}
Let $\varphi\in C_0^\infty(-\tfrac{t}{2},\tfrac{t}{2})$ be radially non-increasing.
Assume $\varphi = const$, $x\in(\tfrac{t}{4}, -\tfrac{t}{4})$, $\int_\mathbb{R} \varphi dx = 1$ and 
define $\varphi^\pm(\cdot) = \varphi(\pm\tfrac{t}{2}+\cdot)$. Then
\begin{equation}\label{eq7.4}
\|\mathcal{T}_{\varphi^\mp,t} \|_{L^p_w(\mathbb{R}_+)\to L^p_w(\mathbb{R}_+)}< \infty,
\end{equation}
provided $1<p<\infty$ and $w\in A_{p,\pm}^{\loc}(\mathbb{R}_+)$.
\end{lemma}
\begin{proof}
(a) Under our assumptions, we have\footnote{This is the one-sided analogue 
of the Proposition in \cite[Section II.2.1]{Stein1993}.}
\begin{equation}\label{eq7.5}
|T_{\varphi^\mp,t} [f]| \le 2\mathcal{M}_{t}^{\pm}[f],
\quad x\in\mathbb{R}_+.
\end{equation}
Indeed, any function $\varphi$ that satisfy the above conditions is uniformly approximated from the above 
by step functions $\varphi_n = \sum_{i=0}^n a_i \chi_{(-t_i, t_i)}$,
where $0<a_i$, $t<4t_i<2t$, and $\int_{\mathbb{R}} \varphi_n dx= 1$. For such functions, we have
\begin{align*}
|\mathcal{T}_{\varphi^\mp,t} [f]|(x) &\le \sum_{i=0}^n a_i \int_{\mp\tfrac{t}{2}-t_i}^{\mp\tfrac{t}{2}+t_i} 
|f|(x - \tau) d\tau\\
&\le \sum_{i=0}^n a_i 2t_i \frac{t+2t_i}{4t_i}\mathcal{M}_{t}^\pm[f](x)
\le 2 \mathcal{M}_{t}^\pm[f](x).
\end{align*}

(b) In view of \eqref{eq7.5} and the inclusion $w\in A_{p,\pm}^{\loc}(\mathbb{R}_+)$, 
the assertion is the direct consequence of \eqref{eq7.3c}.
\end{proof}

Operators $\mathcal{J}_\mp^s$ are invertible for all $s\ge 0$ in the class of smooth functions
restricted to $\mathbb{R}_{\pm}$, respectively (see \cite{Samko}). We denote the
associated inverses by $\bar{\mathcal{J}}_\mp^{-s}$. 
For $0<\varepsilon<1$ and $\varphi^\mp$ from Lemma~\ref{lm7.2}, let 
$\varphi^\mp_\varepsilon(\cdot) = \varepsilon^{-1}\varphi^\mp(\varepsilon^{-1}\cdot)$, 
$\mathcal{J}_{\varepsilon,\mp}^{-s} 
= \bar{\mathcal{J}}_{\mp}^{-s} \mathcal{T}_{\varphi_\varepsilon^\mp,t}$ and
\[
\mathcal{J}_{\mp}^{-s}[f] = \lim_{\varepsilon\to 0} \mathcal{J}_{\varepsilon,\mp}^{-s}[f].
\]
\begin{lemma}\label{lm7.3}
Operators $\mathcal{J}_{\mp}^{-s}:L^{p,s}_w(\mathbb{R}_\pm)\to L^p_w(\mathbb{R}_\pm)$, 
$1< p < \infty$, $s\ge 0$, are one-to-one, provided $w\in A_{p,\pm}^{\loc}(\mathbb{R}_\pm)$.
\end{lemma}
\begin{proof} Straightforward calculations show that $\mathcal{T}_{\varphi_\varepsilon^\mp,t}$ and 
$\mathcal{J}_\mp^s$ commute.\footnote{Note
$\supp \mathcal{F}[\hat{\kappa}^\mp_{s,\ell}], \supp \varphi^\mp \subset\mathbb{R}_\mp$.}
Therefore, for each $f \in L^{p,s}_w(\mathbb{R}_\pm)$  
(by definition $f= \mathcal{J}_\mp^s[\phi]$ for some $\phi \in L^{p}_w(\mathbb{R}_\pm)$), we have
\[
\|\mathcal{J}_{\varepsilon, \mp}^{-s}[f] - \phi\|_{L^p_w(\mathbb{R}_\pm)} 
= \|T_{\varphi_\varepsilon^\mp,t}[\phi] - \phi\|_{L^p_w(\mathbb{R}_\pm)}\to 0,
\quad\text{as}\quad \varepsilon \to 0.
\] 
The conclusion follows from the uniqueness of strong limits.
\end{proof}

Lemma~\ref{lm7.3} indicates that $\mathcal{J}_\mp^s$, $s\ge 0$, are isomorphisms of the scales 
$L^p_w(\mathbb{R}_\pm)$ and $L^{p,s}_w(\mathbb{R}_\pm)$, $1< p< \infty$, 
$w\in A_{p,\pm}^{\loc}(\mathbb{R}_\pm)$. Hence, $L^{p,s}_w(\mathbb{R}_\pm)$, $1< p< \infty$, 
equipped with the norms $\|\cdot\|_{L^{p,s}_w(\mathbb{R}_\pm)} 
:= \|\mathcal{J}_\mp^{-s}[\cdot]\|_{L^{p}_w(\mathbb{R}_\pm)}$, are Banach spaces. 

\subsection{Interpolation} 
Consider a regular vector valued one-sided singular integrals of the form
\begin{subequations}\label{eq7.6}
\begin{equation}\label{eq7.6a}
\mathcal{T}_{\kappa^\pm}[f](x) = (\kappa^{\pm}\ast f)(x),\quad x\in \mathbb{R}_+,
\end{equation}
where $B_0$, $B_1$ are two given Banach spaces, $f:\mathbb{R}_+\to B_0$ and 
$\kappa^\pm(x)\in \mathcal{L}(B_0,B_1)$.
As in the classical theory (see \cite{Stein1993, Graf2009}), we assume 
\begin{equation}\label{eq7.6b}
\mathcal{T}_{\kappa^\pm}\in\mathcal{L}(L^2(\mathbb{R};B_0), L^2(\mathbb{R};B_1)).
\end{equation}
In view of our applications, we consider only compactly supported kernels, 
i.e. kernels with $\supp\kappa^\pm \subset (-t,t)\cap\mathbb{R}_{\pm}$, which 
for all $x,y,\bar{y}\in \supp\kappa^\pm$, with $|x|>0$ and $|y-\bar{y}|\le\tfrac{1}{2}|x-y|$, 
satisfy
\begin{align}
\label{eq7.6c}
&\|\kappa^\pm(x)\|_{B_0\to B_1} \le \tfrac{c}{|x|},\\
\label{eq7.6d}
&\|\kappa^\pm(x-y)-\kappa^\pm(x-\bar{y})\|_{B_0\to B_1} \le 
c\tfrac{|y-\bar{y}|}{|x-y|^2}.
\end{align}
\end{subequations}
In connection with $\mathcal{T}_{\kappa^\pm}$, we define
\begin{equation}\label{eq7.7}
 \mathcal{M}_{\kappa^\pm}[\cdot] = \sup_{\varepsilon>0}
\|\mathcal{T}_{\chi_{|x|>\varepsilon}\kappa^{\pm}}[\cdot]  \|_{B_1}.
\end{equation}
The following result is a straightforward adaptation of the classical "good-$\lambda$ inequality" 
to the one-sided settings, see e.g. \cite[Proposition 6, Section V.4.4]{Stein1993} or \cite[Theorem 9.4.3]{Graf2009}.
\begin{lemma}\label{lm7.4}
Assume $f\in L^1(\mathbb{R}_+)$ satisfy $\supp f \subset \cup_{j} I_j$, where $|I_j|<t$ and
$\dist(I_j, I_k)\ge 2t$, $j\ne k$.
For $\kappa^\pm$ as above and $w\in A^{\loc}_{\infty}$ (see \cite{Ry2001}), 
there exists $0<\alpha_w<1$ such that for any $0<\beta<1$ one can find $\gamma>0$ so that the following holds 
\begin{equation}\label{eq7.8}
w(\{\mathcal{M}_{\kappa^\pm}[f]>\xi\}\cap
\{\mathcal{M}^\mp_{4t}[f] \le \gamma \xi\}) 
\le \alpha_w w(\{\mathcal{M}_{\kappa^\pm}[f]>\beta \xi\}),
\end{equation}
for all $\xi>0$.
\end{lemma}
\begin{proof}
(a) We consider the right-sided operators $\mathcal{M}_{\kappa^-}$, $\mathcal{M}^+_{t}$ only. 
The proof in the left-sided case is identical. Standard arguments (see \cite{Stein1993, Graf2009}) indicate that under 
our assumptions, the level set $E_{\beta \xi}(f) = \{\mathcal{M}_{\kappa^-}[f]>\beta \xi\}$ is open.
The support assumption guarantee that every open connected component $I$ of $E_{\beta \xi}(f)$ 
satisfies $|I|<2t$. It is sufficient to establish \eqref{eq7.8} for a single component $I=(a,b)$, the general 
result follows by summation.

(b) The set $\hat{I} = I / \{\mathcal{M}^+_{4t}[f] > \gamma \xi\})$ is closed in the relative topology
of $I$. If the Lebesgue measure $|\hat{I}|$ of $\hat{I}$ is zero, \eqref{eq7.8} holds trivially. So assume $|\hat{I}|>0$, 
let $x = \min \overline{\hat{I}}$, $\hat{x} = b + (b-x)$, $f_1 = \chi_{[x,\hat{x}]}f$ and $f_2 = (1-\chi_{[x,\hat{x}]})f$
and observe that
\[
w(E_{\xi}(f)\cap I) \le w(E_{\tau \xi}(f_1)\cap I) + w(E_{(1-\tau)\xi}(f_2)\cap I),
\quad 0<\tau<1.
\]
We estimate each term separately.

To bound $w(E_{\tau\xi}(f_1)\cap I)$, we employ the standard weak-type inequality  
(see e.g. \cite[Corollary 2, Section I.7.3]{Stein1993}) to 
obtain initially
\[
\lambda(E_{\tau\xi}(f_1)\cap I) \le \tfrac{c}{\tau\xi} \int_x^{\hat{x}} \| f \|_{B_0} dy
\le \tfrac{2c}{\tau\xi}|I|\mathcal{M}^+_{4t}[f](x) \le \tfrac{2c\gamma}{\tau}|I|,
\]
and then, using the inclusion $w\in A^{\loc}_{\infty}$,
\[
w(E_{\tau\xi}(f_1)\cap I) \le \alpha_w w(I),
\]
with $0<\alpha_w<1$, provided $1<\gamma<0$ is sufficiently small.

(c) To bound $w(E_{\tau\xi}(f_1)\cap I)$, we note that in view of \eqref{eq7.6c} and \eqref{eq7.6d}, 
for $y\in (x,b)$,  we have  
\[
\bigl\|\mathcal{T}_{\chi_{|x|>\varepsilon}\kappa^{-}}[f_2](y) - 
\mathcal{T}_{\chi_{|x|>\varepsilon}\kappa^{-}}[f_2](b)\bigr\|_{B_1} = 0, 
\]
if $\hat{x}\ge b+t$, or 
\begin{align*}
&\bigl\|\mathcal{T}_{\chi_{|x|>\varepsilon}\kappa^{-}}[f_2](y) - 
\mathcal{T}_{\chi_{|x|>\varepsilon}\kappa^{-}}[f_2](b)\bigr\|_{B_1}\\ 
&\qquad \le \int_{\hat{x}+\varepsilon}^{b+t} \|f(z)\|_{B_0} \|\kappa^-(y-z) - \kappa^-(b-z)\|_{B_0\to B_1}dz\\
&\qquad \le c|b-y|\int_{\hat{x}}^{b+t} \|f(z)\|_{B_0} \tfrac{dz}{(z-y)^2}\\
&\qquad \le c \sum_{j\ge 0} \tfrac{|b-y|}{|\hat{x} - y+(2^j-1)(b-y)|^2}
\int_{\hat{x}+(2^j-1) (b-y)}^{\hat{x}+(2^{j+1}-1)(b-y)} \|f\|_{B_0} 
\chi_{[\hat{x},b+t]}dz\\
&\qquad \le 4c \sum_{j\ge 0} \tfrac{ 2^{-j}}{|\hat{x} - x+(2^{j+1}-1)(b-y)|}
\int_{x}^{\hat{x}+(2^{j+1}-1)(b-y)} \|f\|_{B_0} \chi_{[x,b+t]}dz\\
&\qquad \le 8c \mathcal{M}_{3t}^+[f](x)\le 8c \mathcal{M}_{4t}^+[f](x) \le 8c\gamma\xi,
\end{align*}
when $\hat{x}<b+t$. In either case, since $b\notin E_{\beta\xi}(f)$, taking supremum over $\varepsilon>0$, we obtain
\[
\mathcal{M}_{\kappa^-}[f_2](y) \le (\beta+8c\gamma)\xi \le (1-\tau)\xi,\quad y\in(x,b),
\]
provided $0<\gamma<1$ is small and $0<\tau<1$ is chosen appropriately. 
Hence, $w(E_{(1-\tau)\xi}(f_2)\cap I) = 0$ and we conclude that \eqref{eq7.8} holds.
\end{proof}

\begin{corollary}\label{lm7.5}
For $\kappa^{\pm}$ as above and $w\in A^{\loc}_{p,\mp}(\mathbb{R}_+)\cap A^{\loc}_\infty$, 
$1<p<\infty$, the following holds
\begin{subequations}\label{eq7.9}
\begin{align}
\label{eq7.9a}
&\bigl\|\mathcal{M}_{\kappa^\pm}\bigr\|_{L^p_w(\mathbb{R}_+;B_0)
\to L^p(\mathbb{R}_+)} <\infty,\\
\label{eq7.9b}
&\bigl\|\mathcal{T}_{\kappa^\pm}\bigr\|_{L^p_w(\mathbb{R}_+; B_0)
\to L^p_w(\mathbb{R}_+; B_1)} < \infty.
\end{align}
\end{subequations}
\end{corollary}
\begin{proof}
(a) Consider $f\in C_0^\infty(\mathbb{R}_+;B_0)$ initially. Without loss of generality, we may assume 
that $f$ satisfies the support condition of Lemma~\ref{lm7.4} (for any function in $\mathbb{R}$
is a sum of at most four functions satisfying this condition). By our assumptions, 
$g_\varepsilon = \mathcal{T}_{\kappa^{\pm}\chi_{|x|>\varepsilon}}[f]$ is compactly supported and 
smooth, with $\|g_\varepsilon\|_{L^\infty(\mathbb{R}, B_1)}$ bounded independently of $\varepsilon>0$. 
Since $w\in L^{1,\loc}(\mathbb{R}_+)$, we conclude that 
$\|\mathcal{M}_{\kappa^\pm}[f]\|_{L^p_w(\mathbb{R}_+)}<\infty$.
Once this fact is established, we make use of Lemma~\ref{lm7.4} and \eqref{eq7.3c}
to obtain 
\[
\bigl\| \mathcal{M}_{\kappa^\pm}[f] \bigr\|_{L^p_w(\mathbb{R}_+)} \le 
c \bigl\| \mathcal{M}^\mp_{4t}[f] \bigr\|_{L^p_w(\mathbb{R}_+)} \le c' \|f\|_{L^p_w(\mathbb{R}_+;B_0)},
\]
for all $f\in C_0^\infty(\mathbb{R}_+;B_0)$. The standard density argument allows one to pass 
from $C_0^\infty(\mathbb{R}_+;B_0)$ to $L^p_w(\mathbb{R}_+;B_0)$. 
Hence, the bound \eqref{eq7.9a} is settled. 

(b) Estimate \eqref{eq7.9b} is the direct consequence of \eqref{eq7.9a}, as 
\[
\|\mathcal{T}_{\kappa^\pm}[f]\|_{B_1} \le \mathcal{M}_{\kappa^\pm}[f] + 
c\|f\|_{B_0},
\]
a.e. in $\mathbb{R}_+$, see e.g. \cite[Section I.7.4]{Stein1993}. 
\end{proof}

To proceed further, we employ the following local reproducing formula of V. Rychkov (see \cite{Ry2001} for the details)
\begin{subequations}\label{eq7.10}
\begin{equation}\label{eq7.10a}
\delta = \sum_{j\ge 0} \varphi_j^{\pm}\ast\psi_j^{\pm},
\end{equation}
where $\varphi_0^{\pm},\psi^{\pm}\in C_0^\infty(\mathbb{R})$, with  
$\supp \varphi_0^{\pm},\supp\psi^{\pm}\subset (-t,t)\cap\mathbb{R}_{\pm}$ for some $t>0$, 
have non vanishing zeroth moment;
$\varphi^\pm(\cdot) = \varphi_0^\pm(\cdot) - 2^{-1}\phi_0^\pm(2^{-1}\cdot)$, 
$\psi^\pm(\cdot) = \psi_0^\pm(\cdot) - 2^{-1}\psi_0^\pm(2^{-1}\cdot)$ and 
$\varphi^\pm_j(\cdot) = 2^j\varphi^\pm(2^j\cdot)$,
$\psi^\pm_j(\cdot) = 2^j\psi^\pm(2^j\cdot)$, $j\ge 1$. Furthermore, 
both $\varphi_0^\pm$ and $\psi_0^\pm$ can be chosen so that 
\begin{equation}\label{eqeq7.10b}
\int_{\mathbb{R}} x^k \varphi^{\pm} dx =  \int_{\mathbb{R}} x^k \psi^{\pm} dx = 0,\quad
0\le k\le L,
\end{equation}
\end{subequations}
for any given positive integer $L>0$ (in the sequel, we employ symbol $\{\varphi\}_m$ to denote the number 
of vanishing moments of a function $\varphi$).

For $\varphi_0^\pm$ as above, with $\{\varphi^\pm\}_m\ge \max\{0,s\}$, $s\in\mathbb{R}$, we define
\[
\mathcal{S}_{\varphi^\pm}^s[f] 
= \Bigl(\sum_{j\ge 0} 2^{2js}|\varphi_j^{\pm}\ast f|^2\Bigr)^{\frac{1}{2}},\quad x\in\mathbb{R}_+.
\]
\begin{theorem}\label{lm7.6}
For $1< p< \infty$, $s\ge 0$ and $w\in A_{p,+}^{\loc}(\mathbb{R}_+)\cap A^{\loc}_\infty$, we have
\begin{equation}\label{eq7.11}
\|f\|_{L^{p,s}_w(\mathbb{R}_+)} \approx \|\mathcal{S}_{\varphi^-}^s[f]\|_{L^p_w(\mathbb{R}_+)},
\end{equation}
where $\approx$ means the bilateral estimate.
\end{theorem}
\begin{proof}
(a) To begin, we show that 
\begin{equation}\label{eq7.12}
\|f\|_{L^{p}_w(\mathbb{R}_+)} \approx \|\mathcal{S}_{\varphi^\pm}^0[f]\|_{L^p_w(\mathbb{R}_+)},
\end{equation} 
provided $w\in A_{p,\mp}^{\loc}(\mathbb{R}_+)\cap A^{\loc}_\infty$.\footnote{
The proof of \eqref{eq7.12} is identical to that of Theorem 1.10 in \cite{Ry2001}, with the exception 
that, instead of \cite[Theorem 2 and its Corollary, Section V.4.2]{Stein1993}, we invoke Corollary~\ref{lm7.5}.}

Define $\kappa^\pm:\mathbb{R}_+\to \ell_2$ by means of the formulas 
$\kappa^\pm(\cdot) = \{\varphi_j^{\pm}(\cdot)\}_{j\ge 0}$. Operator 
$\mathcal{T}^\pm_{\kappa^\pm}:L^p_w(\mathbb{R}_+)\to L^p_w(\mathbb{R}_+;\ell_2)$ 
fells in the scope of Corollary~\ref{lm7.5}. Hence,
\begin{align*}
\|\mathcal{S}_{\varphi^\pm}^0[f]\|_{L^p_w(\mathbb{R}_+)} 
= \|\mathcal{T}_{\kappa^\pm}f\|_{L^p_w(\mathbb{R}_+;\ell_2)}\le c \|f\|_{L^p_w(\mathbb{R}_+)},
\end{align*}
provided $w\in A_{p,\mp}^{\loc}(\mathbb{R}_+)\cap A^{\loc}_\infty$.

The converse inequality follows from the standard duality argument and the local reproducing formula 
\eqref{eq7.10a}. Indeed, for $g\in L^{p'}(\bar{w}_p)$ supported in $\mathbb{R}_+$, 
we let $g^\pm(\cdot) = g(\pm\cdot)$ and note that $(g^\pm)^- = g^\mp$. Then
\begin{align*}
|\langle f, g\rangle| &= |f\ast g^-|(0) = 
\Bigl|\sum_{j\ge 0} (\varphi_j^\pm\ast f) \ast (\psi_j^\pm\ast g^-)\Bigr|(0)
\le \int_{\mathbb{R}_+}\mathcal{S}_{\varphi^\pm}^0[f] \mathcal{S}_{\psi^\mp}^0[g] d\tau\\
&\le \|\mathcal{S}_{\varphi^\pm}^0[f]\|_{L^p_w(\mathbb{R}_+)} 
\|\mathcal{S}_{\psi^\mp}^0[g]\|_{L^{p'}_{\bar{w}_p}(\mathbb{R}_+)}
\le c \|\mathcal{S}_{\varphi^\pm}^0[f]\|_{L^p(w)}  \|g\|_{L^{p'}_{\bar{w}_p}(\mathbb{R}_+)}.
\end{align*}
Hence, \eqref{eq7.12} is settled.

(b) The general result follows from \cite[Theorem~{2.18}]{Ry2001}, invertibility of $\mathcal{J}_-^\beta$, 
Lem\-ma~\ref{lm7.3} and \eqref{eq7.12}. Indeed, letting
\[
\bar{L}^{p,s}_w(\mathbb{R}_+) = \{ f\,|\,\supp f\in \mathbb{R}_+,\quad 
\| \mathcal{S}_{\varphi^-}^s \|_{L^p_w(\mathbb{R}_+)}<\infty\},
\]
from \cite[Theorem~{2.18}]{Ry2001} we infer that
$\mathcal{J}_-^s\bigl[ \bar{L}^{p,0}_w(\mathbb{R}_+)] 
\subset \bar{L}^{p,s}_w(\mathbb{R}_+)$, $\beta\ge 0$ and the map is onto, as 
$\mathcal{J}_-^s$ is invertible. Hence, by \eqref{eq7.12} and Lemma~\ref{lm7.3},
$\bar{L}^{p,s}_w(\mathbb{R}_+)= \mathcal{J}_-^\beta[\bar{L}^{p,0}_w(\mathbb{R}_+)] 
= \mathcal{J}_-^s[L^{p}_w(\mathbb{R}_+)] =  L^{p,s}_w(\mathbb{R}_+)$ as Banach spaces.
\end{proof}

\begin{remark}\label{lm7.7}
The key feature of Theorem~\ref{lm7.6}, as compared to its analogue \cite[Theorem~{1.10}]{Ry2001}, 
is the use of condition $w\in A_{p,+}^{\loc}(\mathbb{R}_+)\cap A^{\loc}_\infty$, instead of 
$w\in A_p^\loc$. 
The former class is significantly larger than the latter one. For instance, 
$w_\alpha(x) = |x|^\alpha \in A_p^{\loc}$ if and only if $-1<\alpha<\tfrac{p}{p'}$, while 
$w_\alpha(x) \in A_{p,+}^{\loc}(\mathbb{R}_+)\cap A^{\loc}_\infty$, for all $\alpha>-1$.
\end{remark}

In view of Theorem~\ref{lm7.6} and Remark~\ref{lm7.7}, the interpolation identity \eqref{eq2.5} is a 
simple consequence of the  following
\begin{corollary}\label{lm7.8}
Assume $1<p<\infty$ and $w_0, w_1\in A^{\loc}_{p,+}(\mathbb{R}_+)\cap A^{\loc}_{\infty}(\mathbb{R}_+)$.
Then
\begin{equation}\label{eq7.13}
[L^{p,s_0}_{w_0}(\mathbb{R}_+),L^{p,s_1}_{w_1}(\mathbb{R}_+)]_\theta 
= L^{p,(1-\theta)s_0 + \theta s_1}_{w_0^{1-\theta}w_1^\theta}(\mathbb{R}_+),
\; s_0,s_1\ge 0, \; \theta\in(0,1),
\end{equation}
where $[\cdot,\cdot]_\theta$ denotes the standard complex interpolation functor of A. Calderon \cite{BeLo1976}.
\end{corollary}
\begin{proof}
Directly from \eqref{eq7.1} and definition of $A_\infty^\loc$ weights in \cite{Ry2001}, 
it follows that $w_0^{1-\theta}w_1^\theta A^{\loc}_{p,+}(\mathbb{R}_+)\cap A^{\loc}_{\infty}(\mathbb{R}_+)$,
while, in view of Theorem~\ref{lm7.6}, $L^{p,s}_w(\mathbb{R}_+)$ is a retracts of 
$L^{p}_w(\mathbb{R}_+;\ell^{s}_2)$.
Hence, combining the arguments of \cite[Theorem~{5.1.2}]{BeLo1976} 
and \cite[Theorem~{5.5.3}]{BeLo1976}, we have the desired result.
\end{proof}

\bibliographystyle{siamplain}
\bibliography{ms}
\end{document}